%
%
%

%
%

\input amstex
\documentstyle{amsppt}
\input pictex
\magnification=\magstep 1
\pagewidth{160truemm}
\pageheight{240truemm}
\NoBlackBoxes

\font\sss=cmss10

\loadbold

\def\demoend{{\unskip\nobreak\hfill\hbox{\,}\nobreak\hfill
     $\square$}\enddemo}

\subheadskip \bigskipamount

\catcode`\@=11

\def\endremark{\par\revert@envir\endremark\smallskip}
\catcode`\@=12

\font\ss=cmss10

\def\0{{\bold 0}}
\def\1{{\bold 1}}
\def\a{\alpha}
\def\B{\Cal B}

\def\bnd{{\mathop{\text{\bf bnd}}}}
\def\card{\mathop{\text{card}}}
\def\D{\Delta}
\def\d{\delta}
\def\E{\bold E}
\def\e{\varepsilon}

\def\F{\Cal F}
\def\f{\varphi}
\def\ft#1{\foldedtext\foldedwidth{130truemm}{#1}}
\def\G{\Gamma}
\def\g{\gamma}
\def\h{\boldsymbol h}
\def\hh{\bold h}
\def\H{\Cal H}

\def\l{\lambda}
\def\la{\langle}
\def\lra{\longleftrightarrow}

\def\mmu{{\boldsymbol\mu}}
\def\o{\omega}
\def\O{\Omega}
\def\ov{\overline}
\def\P{\bold P}
\def\p{{\frak p}}
\def\part{\partial}
\def\PP{\Cal P}
\def\precc{\preccurlyeq}
\def\Prob{\text{\ss Prob}}
\def\Q{\bold Q}
\def\R{\Bbb R}
\def\ra{\rangle}
\def\RW{\text{\rm RW}}
\def\RWRTP{\text{\rm RWRTP}}
\def\RWTDI{\text{\rm RWTDI}}
\def\s{\sigma}
\def\sd{\rightthreetimes}

\def\supp{{\mathop{\text{\rm supp}}}}

\def\tail{{\mathop{\text{\bf tail}}}}
\def\th{\theta}
\def\Th{\Theta}
\def\toitself{\hookleftarrow}
\def\toto{\mathop{\;\longrightarrow\;}}
\def\wh{\widehat}
\def\wt{\widetilde}
\def\x{\boldsymbol x}

\def\Z{{\Bbb Z}}

\topmatter

\title
Boundaries and harmonic functions for random walks with random transition
probabilities
\endtitle

\author
Vadim A. Kaimanovich$^1$, Yuri Kifer$^2$, Ben-Zion Rubshtein$^3$
\endauthor

\affil
$^1$
CNRS UMR 6625 (IRMAR), Universit\'e Rennes-1, France
\\ \\
$^2$ Institute of Mathematics, Hebrew~University, Jerusalem, Israel
\\ \\
$^3$ School of Mathematics, Ben-Gurion~University, Beer Sheva, Israel
\endaffil

\address
CNRS UMR 6625, IRMAR, Universit\'e Rennes-1, Campus Beaulieu, 35042 Rennes,
France
\endaddress

\email
kaimanov@univ-rennes1.fr
\endemail

\address
Institute of Mathematics, Hebrew University, Jerusalem 96261, Israel
\endaddress

\email
kifer@math.huji.ac.il,
\endemail

\address
School of Mathematics, Ben Gurion University, Beer Sheva, Israel
\endaddress

\email
benzion@cs.bgu.ac.il
\endemail



\dedicatory Dedicated to Eugene Borisovich Dynkin on the occasion of his 80th
anniversary \enddedicatory

\rightheadtext{random walks with random transition probabilities}

\leftheadtext{Vadim A. Kaimanovich, Yuri Kifer, Ben-Zion Rubshtein}

\abstract The usual random walk on a group (homogeneous both in time and in
space) is determined by a probability measure on the group. In a random walk
with random transition probabilities this single measure is replaced with a
stationary sequence of measures, so that the resulting (random) Markov chains
are still space homogeneous, but no longer time homogeneous. We study various
notions of measure theoretical boundaries associated with this model and
establish an analogue of the Poisson formula for (random) bounded harmonic
functions. Under natural conditions on transition probabilities we identify
these boundaries for several classes of groups with hyperbolic properties and
prove the boundary triviality (i.e., the absence of non-constant random bounded
harmonic functions) for groups of subexponential growth, in particular, for
nilpotent groups.
\endabstract

\endtopmatter

\document

\head
Introduction
\endhead

\smallskip

Random walks on groups were intensively studied during the last 40 years (see,
for instance, [Ka96] and the references therein). Their importance is due to
numerous applications, in particular, to the description of boundaries and
spaces of harmonic functions and to the study of ergodic properties of group
actions. Such random walks are Markov chains which are homogeneous both in time
and space and can also be represented as products of independent identically
distributed (i.i.d.) group elements. In the important special case of products
of random matrices additional tools such as Lyapunov exponents can be employed.

A {\it random walk\/} on a group $G$ is determined by a Markov operator
$P=P(\mu)$ which consists in the (right) convolution with a fixed probability
measure $\mu$ on $G$, so that the operator $P$ is invariant with respect to the
action of the group on itself by left translations. There are two models for
further ``randomization'' of these ``ordinary'' random walks. The first model
is usually referred to as {\it random walks in random environment\/} (RWRE) and
consists in considering a probability measure $\l$ on the space of all Markov
operators on $G$. In this model the individual operators (environments) are not
group invariant, although the group structure is taken into account by
requiring the measure $\l$ to be quasi-invariant with respect to the action of
$G$ on the space of environments (more specifically, $\l$ is usually assumed to
be either translation invariant or stationary with respect to the ``moving
environment'' chain, see, for instance, \cite{Kal81}, \cite{KMo84},
\cite{KSi00}). One chooses a random environment according to the distribution
$\l$, and then runs a time (but not space!) homogeneous Markov chain in this
environment.

The other model which we call {\it random walks with random transition
probabilities\/} (RWRTP) is opposite to RWRE in the sense that here one keeps
the space homogeneity but does not assume the time homogeneity. Namely, the
additional randomness is introduced in this model by taking a random sequence
$\mu_0,\mu_1,\dots$ of probability measures on $G$ so that the ($G$-invariant!)
transition probabilities of the arising chain on $G$ at time $n$ are given by
the measure $\mu_n$. The formal description of this model consists in fixing an
invertible ergodic transformation $T$ of a probability space $(\O,\l)$ and a
measurable map $\o\mapsto\mu^\o$. One chooses $\o\in\O$ according to the
distribution $\l$ and then runs the arising {\it random walk with time
dependent increments\/} $\RWTDI(\o)$ determined by the sequence
$$
\mu^\o,\mu^{T\o},\mu^{T^2\o},\dots \;.
$$
The transformation $T$ is usually assumed to be measure preserving, so that the
above random sequence of measures is stationary. In the same way one can also
talk about random sequences of Markov operators on a general state space (which
does not have to be a group or to be endowed with any additional spatial
structure).

This model was first introduced more than 20 years ago in connection with a
model of random automata, and various properties of such Markov chains were
investigated since then in a number of papers (see, for instance, \cite{Or91},
\cite{Ki96} and the references therein). Random walks with random transition
probabilities first appeared in \cite{MR88} (see also \cite{MR94},
\cite{LRW94}, \cite{Ru95}), and products of independent random matrices with
stationarily changing distributions were studied in \cite{Ki01}. Ideologically
and methodically this topic is rather close to {\it random dynamical systems\/}
which were intensively studied in recent years, and Markov chains with random
transition probabilities have the same relation to the classical Markov chains
as random dynamical systems to deterministic ones. In both cases the guiding
philosophy suggests that we have good chances to obtain an additional
non-trivial information about the system if it acquires nice properties after
conditioning by some ergodic stationary process (which we do not have much
information about). Note that in the framework of random walks on groups one
can also make one more step and to combine the RWRE and RWRTP models (so that
individual random chains will be neither space nor time homogeneous).

RWRTP can be considered as a generalization of yet another model of
``randomization'' of the ordinary random walks called {\it random walks with
internal degrees of freedom\/} (RWIDF) \cite{KSz83} or {\it covering Markov
chains\/} \cite{Ka95}. These are $G$-invariant Markov chains on the product of
the group $G$ by another space $X$. The transition probabilities of RWIDF are
$$
p\bigl( (g,x), (gh,y) \bigr) = \ov p(x,y)\mu^{x,y} (h)
$$
(assuming that $X$ is countable), where $\mu^{x,y}$ are probability measures on
$X$, and $\ov p(x,y)$ are the transition probabilities of the quotient chain on
$X$. If the quotient chain has a finite stationary measure, then the associated
RWRTP is determined by the space $\O=X^{\Z}$ endowed with the corresponding
shift-invariant Markov measure and the map
$(\dots,x_{-1},x_0,x_1,\dots)\mapsto\mu^{x_0,x_1}$.

\bigskip

The setup of RWRTP yields a natural notion of {\it random harmonic functions\/}
$f_\o$ on $G$ which satisfy the relation
$$
f_\o(g)=\int f_{T\o}(gh)d\mu^\o(h) \;.
$$
These functions can be considered as harmonic functions of the {\it global\/}
time homogeneous $G$-invariant Markov chain on the product $\O\times G$ with
the transition probabilities
$$
p\bigl((\o,g),(T\o,gh)\bigr)=\mu^\o(h) \;.
$$
This global chain is an immediate analogue of the usual ``space-time'' chain
(the role of ``time'' is played here by the space $\O$ endowed with the
transformation $T$). Therefore, description of all bounded random harmonic
functions amounts to describing the {\it Poisson boundary\/} $\G$ of the global
chain.

In this paper we consider discrete groups $G$ only. We introduce the notion of
the relative (or fiber) Avez type {\it entropy\/} of RWRTP which is close to
the notion of the relative (fiber) entropy in the ergodic theory of random
dynamical systems (this theory is also known under the name of relative ergodic
theory, see \cite{Ki86}). Similarly to the theory of ordinary random walks (see
\cite{KV83}, \cite{Ka00}) we give an entropy criterion for triviality of the
{\it tail boundary\/} of almost all $\RWTDI(\o),\,\o\in\O$, which implies
triviality of the Poisson boundary of the global chain, i.e., absence of
non-trivial bounded random harmonic functions. As a corollary, we prove
convergence of random convolutions to  left invariance for nilpotent groups
(earlier it was established for compact and abelian groups in the works of
Mindlin and Rubshtein \cite{MR94} and of Lin, Rubshtein and Wittman
\cite{LRW94} by completely different methods).

The relationship between the Poisson and the tail boundaries for RWRTP turns
out to be more complicated than for ordinary random walks (where the tail and
the Poisson boundaries coincide with respect to any single point initial
distribution). Indeed, in the RWRTP setup the Poisson boundary does not make
sense for individual $\RWTDI(\o)$ on $G$ (because they are not time
homogeneous). As for the global chain on $\O\times G$, its projection onto $\O$
is deterministic, so that the tail boundary $E$ of the global chain admits a
natural projection onto $\O$ whose fibers are the tail boundaries of
$\RWTDI(\o)$ (in particular, triviality of the tail boundary of almost all
$\RWTDI(\o)$ is equivalent to coincidence of $E$ and $\O$). On the other hand,
the Poisson boundary of any time homogeneous Markov chain is a quotient of its
tail boundary. Therefore, there are two natural projections of the tail
boundary $E$ of the global chain: onto $\O$ and onto the Poisson boundary $\G$.
We say that RWRTP is {\it stable\/} if these two projections separate points of
$E$. If RWRTP is stable, then the tail boundaries of individual $\RWTDI(\o)$
can be identified with the Poisson boundary of the global chain, so that
stability of RWRTP is a property analogous to coincidence of the tail and the
Poisson boundaries for ordinary random walks.

We do not know whether RWRTP on groups are always stable. However, in the final
section under the finite first moment condition we (by using the entropy
technique) explicitly identify the Poisson and the tail boundaries of RWRTP on
discrete groups of isometries of non-positively curved spaces with natural
geometric boundaries (for example, for a free group this natural boundary is
the space of ends). Therefore, these RWRTP are stable.

\bigskip

We do not study here random boundaries for continuous groups which, we
hope, will be dealt with in another paper.

\bigskip

At the end of the paper we put an Appendix devoted to several fundamental
definitions and facts on Borel and Lebesgue spaces, conditional measures,
discrete equivalence relations and ergodic decompositions, which are heavily
used throughout the paper. Although this language has become standard in the
ergodic theory, it may be less known to probabilists, which is why we preferred
to expand on this rather than to restrict ourselves just to an assortment of
references.

\bigskip

{\bf Acknowledgment.} A part of this work was done during the authors'
participation in the ``Random Walks 2001'' program at the Schr\"odinger
International Institute for Theoretical Physics (ESI) in Vienna, whose support
is gratefully acknowledged.

\bigskip

\head
1. Measure theoretical boundaries of Markov chains
\endhead

In this Section we introduce the necessary notations and background from the
general theory of Markov chains, see \cite{Re84}, \cite{Ka92}, \cite{Ki01}.

\subhead
1.1. Markov operators
\endsubhead

\definition{Definition 1.1}
Let $(X,m)$ be a Lebesgue measure space with a $\s$-finite positive measure
$m$. A linear operator $P:L^\infty(X,m)\toitself$ is called
{\it Markov\/} if
\roster
\item"(i)"
$P$ preserves positivity, i.e., $Pf\ge\0$ for any function $f\ge\0$;
\item"(ii)"
$P$ preserves constants, i.e., $P\1=\1$ for the function $\1(x)\equiv 1$;
\item"(iii)"
$P$ is continuous in the sense that $Pf_n\downarrow\0$ a.e. whenever
$f_n\downarrow\0$ a.e.
\endroster
\enddefinition

The adjoint operator $P^\ast$ of a Markov operator $P:L^\infty(X,m)\toitself$
acts on the space of integrable functions on the space $(X,m)$, or, in other
words, on the space of measures $\th$ on $X$ absolutely continuous with respect
to $m$ (notation: $\th\prec m$). We shall use the notation $\th P$ for the
measure on $X$ with the density $P^\ast(d\th/dm)$, so that
$\la\th P,f\ra_m = \la\th,Pf\ra_m$ for any function $f\in L^\infty(X,m)$.

A ($\s$-finite) initial distribution $\th\prec m$ gives rise in a standard
way to a Markov measure $\P_\th$ in the {\it path space\/}
$X^{\Z_+} = \{\x=(x_0,x_1,\dots)\}$ of the {\it associated Markov chain\/} on
$X$. The one-dimensional distributions of $\P_\th$ are $\th P^n$, and the
{\it time shift\/} $(S\x)_n=x_{n+1}$ acts on it as $S(\P_\th)=\P_{\th P}$. A
measure $\th\prec m$ is called a {\it stationary\/} measure of the Markov
operator $P$ if $\th P=\th$, or, equivalently, if the measure $\P_\th$ is
$S$-invariant.

By definition, the conditional expectations $\E_\th$ of the measure $\P_\th$
satisfy the relation $\E_\th \bigl[ f(x_{n+1})|x_n=x \bigr] = Pf(x)$ for any
$n\ge 0$. Since the space $(X,m)$, and therefore all spaces
$(X^\Z_+,\P_\th),\;\th\prec m$ are {\it Lebesgue\/} (see Appendix), these
conditional expectations can be replaced with the integrals with respect to the
corresponding {\it conditional measures\/} $\pi_x$ which are called {\it
one-step transition probabilities\/}. Then the operator $P$ and its adjoint
operator take the form
$$
Pf(x) = \int f(y)\,d\pi_x(y) \;, \qquad
\th P = \int \pi_x\,d\th(x) \;.
\tag 1.1
$$

\remark{Remark 1.2} The measures $\pi_x$ are not necessarily absolutely
continuous with respect to $m$. Still, for any function $f\in L^\infty(X,m)$
the integrals above make sense for $m$-a.e. $x\in X$ by Rokhlin's theorem on
conditional decomposition of measures in Lebesgue spaces (Appendix, Theorem
A.2). We shall use this theorem on several occasions below without further
notice.
\endremark

\subhead
1.2. The tail boundary and harmonic sequences
\endsubhead

Denote by $\a_k,\;k\in\Z_+$ the $k$-th coordinate partition of the path space
$(X^{\Z_+},\P_m)$, and for $0\le k<l\le\infty$ put
$\a_{k,l} = \bigvee_{i=k}^l\a_i$,
i.e., two paths $\x$ and $\x'$ are $\a_{k,l}$-equivalent iff $x_i=x'_i$ for all
$k\le i\le l$. The partitions $\a_{k,\infty},\;k>0$ coincide with the preimage
partitions of the powers $S^k$ of the time shift $S$.

Recall that the measurable partitions of the same space are ordered in such a
way that ``the bigger are the elements, the smaller is the partition''; this
order is denoted by $\precc$. Obviously, $\a_{k+1,\infty}\precc\a_{k,\infty}$
for any $k\in\Z_+$. Let $\a_\infty = \bigwedge_k \a_{k,\infty}$ be the {\it
measurable intersection\/} of the sequence $\a_{k,\infty}$, i.e., the biggest
measurable partition of the space $(X^{\Z_+},\P_m)$ which is smaller than any
partition $\a_{k,\infty}$ (see Appendix for the definition of a measurable
partition). The partition $\a_\infty$ is called the {\it tail partition\/} of
the path space.

\definition{Definition 1.3}
The quotient $E$ of the path space $(X^{\Z_+},\P_m)$ with respect to
the tail partition $\a_\infty$ is called the {\it tail boundary\/}.
Denote by $\tail:X^{\Z_+}\to E$ the corresponding projection.
\enddefinition

\remark{Remark 1.4}
If the transition probabilities $\pi_x$ are a.e. purely atomic (in particular,
if the space $X$ is countable), then the tail boundary can also be described in
terms of the theory of {\it discrete equivalence relations\/} in Lebesgue
spaces (see Appendix for the corresponding definitions). Namely, in this
situation the elements ($\equiv$ equivalence classes) of the partitions
($\equiv$ equivalence relations) $\a_{n,\infty}$ are a.e. countable. Moreover,
the equivalence relations $\a_{n,\infty}$ are {\it discrete\/} in the measure
space $(X^{\Z_+},\P_m)$. The {\it tail equivalence relation\/} is the union
$R=\bigcup_{n\ge 0} \a_{n,\infty}$, so that two sample paths $\x =
(x_0,x_1,\dots)$ and $\x' = (x'_0,x'_1,\dots)$ are $R$-equivalent iff there
exists a number $N\ge 0$ such that $x_n=x'_n$ for all $n\ge N$. The tail
equivalence relation is also discrete, and the tail boundary $E$ is then the
{\it space of the ergodic components\/} of $R$ (see Appendix for the
definition).
\endremark

The space $E$ is endowed with the {\it tail measure type\/} $[\e_m]$ which is
the image of the type of the measure $\P_m$. For any probability measure
$\th\prec m$ the {\it tail measure\/} $\e_\th=\tail(\P_\th)$ is absolutely
continuous with respect to $[\e_m]$. We emphasize that the space $E$ and the
projection $\tail$ are defined in the measure theoretical category, so that
they make sense ``$\P_m$-mod 0'' (i.e., up to the sets of $\P_m$-measure $0$)
only (see Appendix for more details).

The quotient of the path space $X^{\Z_+}$ with respect to the partition
$\a_{n,\infty}$ is the space $X^{[n,\infty)}$ of paths on $X$ running from the
time $n$ only. Therefore, one can consider the space $E$ as the inductive limit
(in the measure theoretical category!) of the sequence of the spaces
$X^{[n,\infty)}$ endowed with the images of the measure $\P_m$. Denote by
$\P_{n,\th}$ the measure on the space $X^{[n,\infty)}$ corresponding to
starting the Markov chain at time $n$ with the initial distribution $\th$.

Projecting the measure $\P_{n,\th}$ onto $E$ gives the associated tail measure
$\e_{n,\th}$. Denote by $\e_{n,x}$ the tail measures on $E$ corresponding to
starting the Markov chain at time $n$ from a point $x\in X$ (cf. Remark 1.2).
Then
$$
\e_{n,x} = \int \e_{n+1,y}\,d\pi_x(y) \;.
\tag 1.2
$$

Since $\tail(\x)=\tail(\x')$ if and only if $\tail(S\x)=\tail(S\x')$, the
action of the time shift $S$ descends from $X^{\Z_+}$ to an invertible
transformation of $E$ (also denoted $S$), and $\e_{n,\th} = S^{-n}\e_\th$.

\definition{Definition 1.5}
A sequence of functions $f_n\in L^\infty(X,m),\;n\in\Z_+$ is called a
{\it harmonic sequence\/} if $f_n=P f_{n+1}$ for any $n\in\Z_+$. Denote by
$HS^\infty(X,m,P)$ the space of harmonic sequences endowed with the norm
$\sup_n \|f_n\|_\infty$.
\enddefinition

\proclaim{Theorem 1.6 {\rm (\cite{Re84}, \cite{Ka92})}}
The spaces $HS^\infty(X,m,P)$ and $L^\infty(E,[\e_m])$ are isometric. This
isometry is established by the formulas
$$
\lim_{n\to\infty} f_n(x_n) = \wh f(\tail(\x)) \;, \qquad
f_n(x) = \la \wh f, \e_{n,x} \ra \;,
\tag 1.3
$$
where $\{f_n\}\in HS^\infty(X,m,P)$ and $\wh f\in L^\infty(E,[\e_m])$.
\endproclaim

\subhead
1.3. The Poisson boundary and harmonic functions
\endsubhead

\definition{Definition 1.7}
The space $\G$ of ergodic components of the shift $S$ in the path space
$(X^{\Z_+},\P_m)$ is called the {\it Poisson boundary\/} of the Markov operator
$P$. Denote the corresponding projection by $\bnd:X^{\Z_+}\to\G$.
\enddefinition

The Poisson boundary can also be defined as the space of ergodic components of
the transformation $S$ (induced by the shift in the path space) of the tail
boundary $E$. Therefore, the map $\bnd$ is the result of the composition of the
maps $\tail:X^{\Z_+}\to E$ and $\p_\G:E\to\G$. Denote by
$[\nu_m]=\bnd([\P_m])=\p_\G([\e_m])$ the {\it harmonic measure type\/} on $\G$,
and by $\nu_\th=\bnd(\P_\th)=\p_\G(\e_\th)$ the {\it harmonic measure\/}
corresponding to an initial probability distribution $\th\prec m$, so that
$\nu_\th\prec[\nu_m]$. By $\nu_x$ we shall denote the harmonic measures
corresponding to individual points $x\in X$ (cf. Remark 1.2). Since $\G$ is the
space of ergodic components of the action of $S$ on $E$, for any $n\ge 0$ we
have $\p_\G(\e_{n,\th}) = \p_\G(\e_\th) = \nu_\th$. Therefore, (1.2) implies
that the harmonic measures on $\G$ satisfy the {\it stationarity relation\/}
$$
\nu_x = \int \nu_y\,d\pi_x(y) \;.
$$

\definition{Definition 1.8}
A function $f\in L^\infty(X,\mu)$ is called {\it harmonic\/} with respect to a
Markov operator $P:L^\infty(X,m)\toitself$ if $f=Pf$. Denote by
$H^\infty(X,m,P)$ the subspace of $L^\infty(X,m)$ consisting
of harmonic functions.
\enddefinition

Any function $f\in H^\infty(X,m,P)$ determines the harmonic sequence $f_n\equiv
f$, so that $H^\infty(X,m,P)$ is isometrically embedded into the space
$HS^\infty(X,m,P)$. Since the subspace of $S$-invariant functions in
$L^\infty(E,[\e_m])$ is naturally isometric to the space $L^\infty(\G,[\nu_m])$
of all bounded measurable functions on the Poisson boundary (which is the space
of $S$-ergodic components in $E$), Theorem 1.6 implies

\proclaim{Theorem 1.9 {\rm (\cite{Re84}, \cite{Ka92}) }}
The spaces $H^\infty(X,m,P)$ and $L^\infty(\G,[\nu_m])$ are isometric.
The isometry is established by the formulas
$$
\lim_{n\to\infty} f_n(x_n) = \wh f(\bnd(\x)) \;, \qquad
f(x) = \la \wh f, \nu_x \ra \;,
\tag 1.4
$$
where $f\in H^\infty(X,m,P)$ and $\wh f\in L^\infty(\G,[\nu_m])$.
\endproclaim

Formula (1.4) and its time dependent counterpart (1.3) are called the
{\it Poisson formulas\/}. See \cite{Fu63}, \cite{Ka96} for a relationship with
the classical Poisson formula for bounded harmonic functions on the unit disk
(which is the origin of this term). Criteria of triviality and of coincidence
of the tail and the Poisson boundaries for general Markov chains are given by
the {\it 0--2 laws\/} \cite{De76}, \cite{Ka92}.

\subhead
1.4. Non-homogeneous Markov chains
\endsubhead

The notions introduced above also apply to Markov chains which are {\it not\/}
homogeneous in time. In this situation instead of a single Markov operator $P$
we have a sequence of Markov operators $P_n:L^\infty(X,m)\toitself$
on the same space $(X,m)$. The operator $P_n$ governs the transitions of the
Markov chain at time $n$, so that the measure $\P_\th$ in the path space
$X^{\Z_+}$ corresponding to an initial distribution $\th$ on $X$ satisfies the
relations
$$
\E_\th ( f(x_{n+1}) | x_n=x ) = P_n f(x) \;.
$$
The one-dimensional distribution of $\P_\th$ at time $n+1$ is $\th P_{0,n}$,
where
$$
P_{k,n} = P_k P_{k+1} \cdots P_n \;,\qquad 0\le k\le n
\tag 1.5
$$
are the ``time $k$ to time $n+1$'' transition operators. The standard way to
``make'' such chains homogeneous consists in extending the state space by
passing to the ``space-time'' $\Z_+\times X$ (or to $\Z\times X$ when dealing
with negative times as well). Then one can talk about a single {\it space-time
operator\/}
$$
Pf(n,\cdot)=P_n f(n+1,\cdot)
\tag 1.6
$$
on $\Z_+\times X$ and about the corresponding time homogeneous Markov chain
(which is called the {\it space-time chain\/}). The projection of the
space-time chain onto $\Z_+$ is deterministic and consists in moving forward
with unit speed.

The notions of the tail boundary and of harmonic sequences carry over to
non-homogeneous Markov chains without any changes, whereas the Poisson boundary
and harmonic functions do not make much sense in this situation. For the
space-time chain the tail boundary coincides with the Poisson boundary and is
the product of the tail boundary of the original non-homogeneous chain by $\Z$,
see \cite{Ka92} for more details.

\subhead
1.5. Random Markov operators
\endsubhead

\definition{Definition 1.10}
A {\it random Markov operator\/} on a space $(X,m)$ is determined by a measure
type preserving transformation $T$ of a probability space $(\O,\l)$ and a
measurable map $\o\mapsto P_\o$ from $\O$ to the space of Markov operators on
$(X,m)$. We shall call $(X,m)$ the {\it state space\/} and $(\O,\l)$ the
{\it base space\/} of the random Markov operator $\{P_\o\}$. Here by
measurability of the map $\o\mapsto P_\o$ we mean that the integral
$\la P_\o f,g\ra_m$ is a measurable function of $\o$ for any two functions
$f\in L^\infty(X,m), g\in L^1(X,m)$.
\enddefinition

$$
\ft{
For simplicity we shall always assume that the transformation $T$ is ergodic
and invertible. In most applications the measure $\l$ on $\O$ is in addition
assumed finite and $T$-invariant in order to guarantee the ``stochastic
homogeneity'' of the sequence of operators $P_{T^n\o}$.
}
$$

\medskip

For any $\o\in\O$ we have a non-homogeneous Markov chain on $X$ determined by
the sequence of operators $P_\o, P_{T\o}, P_{T^2\o}, \dots$. Denote by
$\P_{\o,\th}$ the measure in its path space $X^{\Z_+}$ corresponding to
an initial distribution $\th$ on $X$, and by $E_\o$ its tail boundary. By
$\e_{\o,\th}=\tail_\o\P_{\o,\th}$ we denote the tail measure on $E_\o$
corresponding to the initial distribution $\th$ on $X$ (the subscript $\o$
indicates that the map $\tail_\o$ is defined on the path space of the chain
determined by $\o$).

Simultaneously with the Markov operators $P_\o:L^\infty(X,m)\toitself$ we shall
also consider the ``global'' Markov operator
$$
Pf(\o,\cdot) = P_\o f(T\o,\cdot) \;.
\tag 1.7
$$
acting on the space $L^\infty(\O\times X,\l\otimes m)$. The operator $P$ is an
immediate analogue of the space-time operator (1.6), the only difference
being that the role of ``time'' here is played by the space $\O$ endowed with
the transformation $T$. The transition probabilities of the operator $P$ are
$$
\pi_{\o,x} = \d_{T\o} \otimes \pi_x(\o) \;,
$$
where $\pi_x(\o)$ are the transition probabilities of the operator $P_\o$. The
sample paths of the operator $P$ have the form
$$
\ov\x = (\ov x_0,\ov x_1,\dots) \;, \qquad \ov x_n = (T^n\o,x_n) \;,
$$
where $\x=(x_n)$ is a sample path of the non-homogeneous Markov chain
determined by the sequence of operators $(P_\o,P_{T\o},\dots)$. Therefore, the
path space of the operator $P$ can be identified with $\O\times X^{\Z_+}$ by
the map
$$
\Pi: \ov\x  \mapsto (\o,\x) \;,\qquad \x=(x_0,x_1,\dots)\in X^{\Z_+} \;.
\tag 1.8
$$
As usually, denote by $\P_{\ov\th}$ the measure on the path space
$(\O\times X)^{\Z_+}$ of the operator $P$ corresponding to an initial
distribution $\ov\th$ on $\O\times X$. Then
$$
\Pi\P_{\ov\th} = \int \d_\o\otimes \P_{\o,\th_\o}\,d\th(\o) \;,
\tag 1.9
$$
where $\th$ is the image of $\ov\th$ under the projection from $\O\times X$
onto $\O$, and $\th_\o,\;\o\in\O$ are the conditional measures of this
projection.

Denote by $E$ the tail boundary of the operator $P$ (1.7).
Let
$$
\Pi_\O:\ov\x\mapsto\o \;,\qquad \ov\x=\bigl((\o,x_0),(T\o,x_1),\dots\bigr)
\tag 1.10
$$
be the composition of the map $\Pi$ and the projection from $\O\times X^{\Z_+}$
onto $\O$. Since the transformation $T$ is invertible, $\Pi_\O$ is measurable
with respect to the tail partition of the path space $(\O\times X)^{\Z_+}$.
Therefore, $\Pi_\O$ determines a natural projection $\p_\O:E\to\O$.

\proclaim{Proposition 1.11}
The fibers of the projection $\p_\O$ are the tail boundaries $E_\o$ of the
non-homogeneous Markov chains on $X$ associated with the points $\o\in\O$. More
precisely, for an arbitrary initial distribution $\ov\th\prec\l\otimes m$ on
$\O\times X$ denote by $\th_\o,\,\o\in\O$ its conditional measures on $X$. Then
the conditional measures of the tail measure $\e_{\ov\th}$ on $E$ with respect
to the projection $\p_\O$ coincide with the tail measures $\e_{\th_\o}$ on the
tail boundaries $E_\o,\,\o\in \O$.
\endproclaim

\demo{Proof}
Although this claim is true in stated full generality, we shall only need it in
the situation when the state space $X$ is countable, so that our proof is
restricted just to this case where it is a direct consequence of the
transitivity of ergodic decompositions of discrete equivalence relations
(Appendix, Theorem A.7) [Actually, the same easy argument is applicable to
general Lebesgue spaces $X$ as well, provided the transition probabilities
$\pi_x(\o)$ are purely atomic.] Namely, we shall use the fact that under this
assumption the tail boundary $E$ of the operator $P$ is the space of the
ergodic components of the tail equivalence relation $R$ in the path space
$\bigl((\O\times X)^{\Z_+},\P_{\ov\th}\bigr)$ (see Remark 1.4), whereas the
tail boundaries $E_\o$ are the ergodic components of the restrictions of $R$ to
the fibers of the projection $\Pi_\O$ (1.10), i.e., to the path spaces of the
operators $P_\o$. Therefore, we are in the setup of Theorem A.7. More
precisely, we have the diagram
$$
\beginpicture
 \setcoordinatesystem units <1mm,1mm>
 \setplotarea x from -4 to 50, y from -20 to 0
 \put{$\bigl((\O\times X)^{\Z_+},\P_{\ov\th}\bigr)$} at 0 0
 \put{$(E,\e_{\ov\th})$} at -1 -16
 \put{$(\O,\th)\qquad,$} at 40 -16
 \put{$\tail$} at -6 -7
 \put{$\Pi_\O$} at 23 -5
 \put{$\p_\O$} at 20 -14
 \arrow <3.6truept> [.4,1.75] from -1 -4 to -1 -12
 \arrow <3.6truept> [.4,1.75] from 14 -4 to 28 -12
 \arrow <3.6truept> [.4,1.75] from 14 -16 to 28 -16
\endpicture
$$
with $\tail$ being the projection of the space $\bigl((\O\times
X)^{\Z_+},\P_{\ov\th}\bigr)$ onto the space $(E,\e_{\ov\th})$ of the ergodic
components of the equivalence relation $R$. Thus, the claim of Proposition~1.11
is a specification of the claim of Theorem A.7.
\demoend

\pagebreak

Denote by $\G$ the Poisson boundary of the operator $P$. By definition, there
is a projection $\p_\G:E\to\G$. Let $\s_\O$ and $\s_\G$ be the preimage
partitions of the tail boundary $E$ determined by the projections $\p_\O$ and
$\p_\G$, respectively.

$$
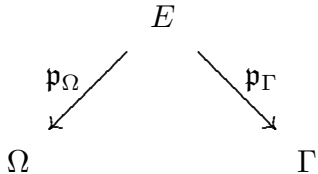
\beginpicture
\setcoordinatesystem units <1mm,1mm>
\setplotarea x from -16 to 16, y from -20 to 0
\put{$E$} at 0 0
\put{$\O$} at -16 -16
\put{$\G$} at 16 -16
\put{$\p_\O$} at -11 -7
\put{$\p_\G$} at 11 -7
\arrow <3.6truept> [.4,1.75] from -4 -4 to -12.5 -12.5
\arrow <3.6truept> [.4,1.75] from 4 -4 to 12.5 -12.5
\endpicture
$$

\definition{Definition 1.12}
A random Markov operator $\{P_\o\}$ has a {\it stable Poisson boundary\/} if
the common refinement $\s_\O\vee\s_\G$ of the partitions $\s_\O$ and $\s_\G$
coincides with the point partition $\s_E$ of the tail boundary $E$, i.e., if
the projections $\p_\O:E\to\O$ and $\p_\G:E\to\G$ separate points of $E$.
\enddefinition

If $\{P_\o\}$ has a stable Poisson boundary $\G$, then the tail boundaries $E$
and $E_\o,\;\o\in\O$ can be identified with the product $\O\times\G$ and with
$\G$, respectively. Therefore, in this situation the same space $\G$ is
responsible (via the Poisson formula) for an integral representation both of
$P$-harmonic functions on $\O\times X$ and of $(P_\o,P_{T\o},\dots)$-harmonic
sequences on $X$ for a.e. $\o\in\O$. In other words, the boundary behaviour of
the operators $P_\o$ and the operator $P$ (in the latter case modulo dependence
on the initial state) is described by the same space $\G$. If the Poisson
boundary $\G$ of the operator $P$ is trivial, then stability of $\{P_\o\}$
means that the tail boundary $E$ of $P$ coincides with $\O$. The general case
can be reduced to this situation by conditioning the operator $P$ by points of
$\G$.

One can easily give a simple (if somewhat degenerate) example of a random
Markov operator whose Poisson boundary is unstable in the sense of Definition
1.12. Basically, it consists just in taking an ergodic skew product over an
ergodic invertible transformation. Indeed, consider on $\O\times X$ the skew
product  transformation $\wt T(\o,x)=(T\o,\f(\o,x))$ with the base $T$, where
$\f:\O\times X\to X$ is a measurable map such that $\f(\o,\cdot):X\to X$ is
invertible, and let $P_\o f(x)=f(\f(\o,x))$ be the associated family of
deterministic Markov operators. Then the operator $P$ is also deterministic and
corresponds to the transformation $\wt T$. Therefore, the tail boundary of $P$
is $E=\O\times X$, whereas the tail boundaries of the operators $P_\o$ are
$E_\o=X$. If $\wt T$ is ergodic (e.g., see \cite{CFS82} for examples), then the
Poisson boundary $\G$ of $P$ is trivial, which gives an example we are looking
for.

It would be interesting to find general sufficient conditions for stability of
the Poisson boundary of random Markov operators. In Section 4 we shall prove it
for group invariant random Markov operators on certain classes of groups.

\head
2. Random walks with time dependent increments
\endhead

In this Section we consider random walks on groups which are {\it space
homogeneous\/} without being {\it time homogeneous\/}.

\pagebreak

\subhead
2.1. Definitions and notations
\endsubhead

Let $G$ be a countable group, and $m=m_G$ be the counting measure on $G$.
Denote by $\PP G$ the space of probability measures on $G$. Any measure
$\mu\in\PP G$ determines the Markov operator
$$
P_\mu f(g) = \sum_h \mu(h) f(gh)
$$
commuting with the action of the group $G$ on itself by left translations. The
operator $P_\mu$ acts on measures on $G$ by (right) convolution with $\mu$:
$$
\th P_\mu=\th\mu \;.
$$
The associated Markov chain on $G$ with the transition probabilities
$$
\Prob(x_{n+1}=gh | x_n=g) = \mu(h)
$$
which is homogeneous both in time and in space is called the (right)
{\it random walk\/} (RW) on $G$ determined by the measure $\mu$ and is  denoted
$\RW(\mu)$ (e.g., see \cite{KV83}).

\definition{Definition 2.1}
For any sequence $\mmu=\{\mu_0,\mu_1,\dots\}\in\PP G^{\Z_+}$ the associated
sequence of Markov operators $P_n=P_{\mu_n}$ determines the Markov chain on
$G$ whose transition probabilities at time $n$ are
$$
\Prob(x_{n+1}=gh | x_n=g)
= \mu_n(h) \;.
$$
This chain is called a {\it random walk with time dependent increments\/}
(RWTDI) on $G$, and we denote it $\RWTDI(\mmu)$.
\enddefinition

Random walks with time dependent increments are {\it homogeneous in space, but
not in time\/} unless the sequence $\mmu$ is constant, in which case we have a
usual random walk on the group $G$ homogeneous both in time and space. All
objects connected with $\RWTDI(\mmu)$ obviously depend on the sequence $\mmu$.
However, for the sake of keeping the notations concise, we shall usually omit
the argument $\mmu$.

The ``time $k$ to time $n+1$'' transition operators (1.5) of $\RWTDI(\mmu)$
are
$$
P_{k,n} f(g)
= P_{\mu_k}\cdots P_{\mu_n} f(g)
= \sum_h \mu_{k,n}(h) f(gh)
= P_{\mu_{k,n}} f(g) \;,
$$
where
$$
\mu_{k,n}=\mu_k\mu_{k+1}\dots\mu_n
$$
is the convolution of the measures $\mu_k,\mu_{k+1},\dots,\mu_n$.

Denote by $\P_{n,\th}$ the measure in the path space $G^{[n,\infty)}$
corresponding to starting $\RWTDI(\mmu)$ at time $n\in\Z_+$ with an initial
distribution $\th$. If $\th=\d_g,\,g\in G$, then we use the notation
$\P_{n,g}$. We shall omit the subscript $n$ if $n=0$ and the subscript $g$ if
$g=e$. In particular, we denote by $\P$ the probability measure on $G^{\Z_+}$
corresponding to starting $\RWTDI(\mmu)$ from the group identity at time 0.
Since $G$ acts on the path space $G^{\Z_+}$ coordinate-wise as $(g\x)_n=gx_n$,
and the transition probabilities of $\RWTDI(\mmu)$ are $G$-invariant,
$\P_{n,g}=g\P_n,\;g\in G$, and $\P_{n,\th}=\th\P_n$ for any initial
distribution $\th$.

\subhead
2.2. The tail boundary and conditional chains
\endsubhead

Denote by $E$ the tail boundary of $\RWTDI(\mmu)$, and by
$\e_{n,\th}=\tail(\P_{n,\th})$ (resp., $\e_{n,g}$, etc.) the tail measures on
$E$. There are two types of the tail behaviour of $\RWTDI(\mmu)$: one (rather
obvious) is connected with the dependence on the starting point and can be
dealt with by passing to a smaller group (see Theorem 2.10 below); the other
one is more interesting and reflects the ``true'' tail behaviour which can not
be reduced to a dependence on the initial state. Because of the space
homogeneity of RWTDI, for the study of the latter we may always assume that the
starting point is the group identity. Denote by $\e=\tail(\P)$ the associated
tail measure. Below by the {\it tail boundary\/} of $\RWTDI(\mmu)$ we shall
always mean the space $(E,\e)$ (sometimes we shall also call it the {\it local
tail boundary\/}), whereas the space $E$ endowed with measure type $[\e_m]$
will be referred to as the {\it total tail boundary\/}.

The action of $G$ on the path space $G^{\Z_+}$ commutes with the time shift, so
that this action descends to an action on $E$ which preserves the tail measure
type $[\e_m]$, and $\e_{n,g\th}=g\e_{n,\th}$ for any $n\in\Z_+,g\in G$ and any
measure $\th$ on $G$. In particular, $\e_{n,\th}=\th\e_n$. The stationarity
relations (1.2) then imply that
$$
\e_n = \sum \mu_n(g) g\e_{n+1} = \mu_n\e_{n+1} \;.
\tag 2.1
$$

\proclaim{Proposition 2.2}
The family of measures $\P^\g,\g\in E$ on $G^{\Z_+}$ defined on cylinder sets
$$
C_{e,g_1,\dots,g_n}=\{\x\in G^{\Z_+}: x_0=e,x_1=g_1,\dots,x_n=g_n\}
$$
as
$$
\P^\g(C_{e,g_1,\dots,g_n})
= \P(C_{e,g_1,\dots,g_n}) \frac{d g_n\e_n}{d\e}(\g)
$$
is the canonical system of conditional measures of the measure $\P$ with
respect to the tail boundary.
\endproclaim

\demo{Proof}
If $A$ is a measurable subset of the tail boundary $E$ with
$\e(A)=\P(\tail^{-1}A)>0$, then by the Markov property
$$
\P (C_{e,g_1,\dots,g_n} \cap \tail^{-1}A)
=\P (C_{e,g_1,\dots,g_n}) \P_{n,g_n}(\tail^{-1}A)
=\P (C_{e,g_1,\dots,g_n}) g_n \e_n(A)
$$
for any cylinder set $C_{e,g_1,\dots,g_n}$, whence for the conditional measure
$\P^A(\cdot)=\P(\cdot|\tail^{-1}A)$ we have
$$
\P^A(C_{e,g_1,\dots,g_n})
= \frac{\P(C_{e,g_1,\dots,g_n})g_n\e_n(A)}{\P(\tail^{-1}A)}
= \P(C_{e,g_1,\dots,g_n})\frac{g_n\e_n(A)}{\e(A)} \;.
$$
Therefore,
$$
\P^A(C_{e,g_1,\dots,g_n})
= \frac1{\e(A)} \int_A \P^\g(C_{e,g_1,\dots,g_n})\,d\e(\g)
$$
for any measurable subset $A\subset E$, which implies the claim.
\demoend

\remark{Remark 2.3}
The definition of conditional measures $\P^\g,\,\g\in E$ can be rewritten as
$$
\align
\P^\g(C_{e,g_1,\dots,g_n})
& = \P(C_{e,g_1,\dots,g_n}) \frac{d g_n\e_n}{d\e}(\g) \\
&= \P(C_{e,g_1,\dots,g_n}) \frac{d g_1\e_1}{d\e}(\g)
\frac{d g_2\e_2}{d g_1\e_1}(\g) \cdots \frac{d g_n\e_n}{d g_{n-1}\e_{n-1}}(\g)
\;.
\endalign
$$
Thus, the measures $\P^\g$ are the measures in the path space of
{\it conditional Markov chains\/} on $G$ with the transition probabilities
$$
\aligned
\Prob^\g (x_{n+1}=g_{n+1} | x_n=g_n)
&= \Prob (x_{n+1}=g_{n+1} | x_n=g_n) \frac{d g_{n+1}\e_{n+1}}{d g_n\e_n}(\g) \\
&= \mu_n(g_n^{-1} g_{n+1}) \frac{d g_{n+1}\e_{n+1}}{d g_n\e_n}(\g) \;.
\endaligned
\tag 2.2
$$
These conditional chains are, generally speaking, inhomogeneous both in space
and in time. Note that (2.2) makes sense only when
$g_{n+1}\e_{n+1}\prec g_n\e_n$. However, as it follows from (2.1), this
relation is satisfied whenever $\mu_{0,n-1}(g_n)$ and $\mu_n(g_n^{-1}g_{n+1})$
are both non-zero, i.e., the conditional chains are well-defined on the whole
{\it attainability space-time cone\/} in $\Z_+\times G$ with the origin $(0,e)$
$$
C = C(0,e) = \{ (n+1,x): n\in\Z_+,\;\mu_{0,n}(x)>0 \} \;.
\tag 2.3
$$
\endremark

Given a measurable partition $\xi$ of the total tail boundary $(E,[\e_m])$
denote by $E^\xi$ the associated quotient space, and by $\e^\xi_{n,\th}$, etc.
the images of the corresponding tail measures under the projection
$\g\mapsto\xi(\g)$ from $E$ to $E^\xi$. A partition $\xi$ is called
{\it $G$-invariant\/} if the action of $G$ on $E$ maps elements of $\xi$ onto
elements of $\xi$ (although individual elements of $\xi$ do not have to be
fixed by the action). If $\xi$ is a $G$-invariant partition, then the action
of $G$ descends from $E$ to the corresponding quotient space $E^\xi$.
Reproducing the proof of Proposition 2.2 we get

\proclaim{Proposition 2.4}
Let $\xi$ be a $G$-invariant measurable partition of the tail boundary $E$.
Then the family of measures $\P^{\xi(\g)},\g\in E$ on $G^{\Z_+}$ defined on
cylinder sets as
$$
\P^{\xi(\g)}(C_{e,g_1,\dots,g_n})
= \P(C_{e,g_1,\dots,g_n}) \frac{d g_n\e^\xi_n}{d\e^\xi}(\xi(\g))
$$
is the canonical system of conditional measures of the measure $\P$ with
respect to the quotient $E^\xi$ of the tail boundary by the partition $\xi$.
\endproclaim

\subhead
2.3. Triviality of the tail boundary
\endsubhead

\definition{Definition 2.5}
If the tail boundary $(E,\e)$ of $\RWTDI(\mmu)$ is a singleton, we shall say
that it is {\it trivial\/}. In this Section we shall also use the term
{\it local triviality\/} in order to distinguish it from the {\it total
triviality\/} of the tail boundary when the total tail boundary $(E,[\e_m])$)
is a singleton.
\enddefinition

\proclaim{Theorem 2.6 {\rm (\cite{De76}, \cite{Ka92}) }}
The tail boundary of $\RWTDI(\mmu)$ is totally trivial iff
$$
\| g\mu_{m,n} - \mu_{0,n} \| \toto_{n\to\infty} 0
\qquad\forall\,g\in G,\;m\in\Z_+ \;,
$$
and it is locally trivial iff
$$
\| g\mu_{m,n} - \mu_{0,n} \| \toto_{n\to\infty} 0
\qquad\forall\,g\in\supp\mu_{0,m-1},\;m\in\Z_+ \;,
$$
where $\|\th\|$ denotes the total variation of a measure $\th$.
\endproclaim

\proclaim{Corollary} If the tail boundary of $\RWTDI(\mmu)$ is totally trivial,
then for any $m\in\Z_+$ the sequence of measures $\mu_{m,n}$ strongly converges
(as $n$ tends to infinity) to left invariance on $G$, and therefore the group
$G$ must be amenable.
\endproclaim

\remark{Remark 2.7}
Theorem 2.6 implies that the total triviality of the tail boundary of
$\RWTDI(\mmu)$ is equivalent to the total triviality of the tail boundary of
$\RWTDI(S\mmu)$, where $S\mmu=(\mu_1,\mu_2,\dots)$ is the shift of
$\mmu=(\mu_0,\mu_1,\dots)$. It also implies that local triviality of the tail
boundary of $\RWTDI(S\mmu)$ follows from local triviality of the tail boundary
of $\RWTDI(\mmu)$. However, the converse is not true in general. For the
simplest example take for $\mu_0$ any measure whose support consists of more
than one point, and put $\mu_1=\mu_2=\dots=\d_e$.
\endremark

Obviously, total triviality implies local triviality. We shall show that
the converse is also true under natural irreducibility conditions. For any
given $g\in G$ the sequence $\| g\mu_{0,n} - \mu_{0,n}\|$ is clearly
non-decreasing. Denote by $\D(g)=\D(g,\mmu)$ its limit.

\proclaim{Proposition 2.8}
If the tail boundary of $\RWTDI(\mmu)$ is locally trivial, then $\D(g)$ equals
either $0$ or $2$ for any $g\in G$.
\endproclaim

\demo{Proof}
Let
$$
A = \{ \x\in G^{\Z_+} : (n,x_n)\in C \;\text{for a certain}\; n\in\Z_+ \} \;,
$$
where $C=C(0,e)$ is the attainability cone (2.3) in $\Z_+\times G$. Clearly,
for $\P_m$-a.e. path $\x$ if $(k,x_k)\in C$ then also $(n,x_n)\in C$ for all
$n>k$, so that the set $A$ is measurable with respect to the tail partition.
Since the tail boundary is locally trivial, $\P_g(A)$ equals 0 or 1 for any
$g\in G$.

Suppose that $\D(g)<2$, i.e., for a certain $n\ge 0$ the measures $g\mu_{0,n}$
and $\mu_{0,n}$ are non-singular. Then $\P_g(A)>0$, and by the above
$\P_g(A)=1$, so that $\P_g$-a.e. path eventually hits the cone $C$. Denote by
$$
\tau(\x) = \min\{n\in\Z_+: (n,x_n)\in C\}
$$
the first hitting time, and by $\th$ the corresponding hitting distribution on
$C$:
$$
\th(n,x) = \P_g \{\x\in G^{\Z_+}: \tau(\x)=n, x_n=x \} \;.
$$
By the Markov property for any $n\in\Z_+$ the measure $g\mu_{0,n}$ decomposes
as
$$
g\mu_{0,n} = \sum_{(k,x)\in C :k\le n} \th(k,x) x\mu_{k,n} + \a_n \;,
$$
where $\|\a_n\|\to 0$. On the other hand, by the local triviality of the tail
boundary
$$
\|x\mu_{k,n} - \mu_{0,n}\| \to 0
$$
for any $(k,x)\in C$, and we are done.
\demoend

Denote by $G(\mmu)$ the subgroup of $G$ generated by all $g\in G$ such that the
measures $\mu_{0,n}$ and $g\mu_{0,n}$ are non-singular for a certain $n$.
Proposition 2.8 implies

\proclaim{Proposition 2.9}
If the tail boundary of $\RWTDI(\mmu)$ is locally trivial, then
$$
G(\mmu) = \{g\in G: \D(g)=0 \} \;.
$$
\endproclaim

\proclaim{Corollary}
If the tail boundary of $\RWTDI(\mmu)$ is locally trivial, then the group
$G(\mmu)$ is amenable.
\endproclaim

\proclaim{Theorem 2.10}
If the tail boundary of $\RWTDI(\mmu)$ is locally trivial, then its total tail
boundary is isomorphic to the coset space $G/G(\mmu)$, and the boundary map
$\tail$ has the form
$$
\tail(\x) = x_0 G(\mmu) \;.
$$
\endproclaim

\demo{Proof}
As it follows from Proposition 2.9, sample paths of $\RWTDI(\mmu)$ issued from
different cosets of $G(\mmu)$ never intersect, so that the map
$\x\mapsto x_0 G(\mmu)$ is indeed measurable with respect to the tail
partition. On the other hand, again by Proposition 2.9, the tail boundary is
trivial with respect to any initial distribution concentrated on a single
coset.
\demoend

We shall say that $\RWTDI(\mmu)$ is {\it irreducible\/} if $G(\mmu)=G$. In
particular, if all points of $G$ are attainable with positive probability
from the group identity ($\equiv$ from an arbitrary starting point), i.e., if
$\bigcup_{n\in\Z_+}\supp\mu_{0,n}=G$, then $\RWTDI(\mmu)$ is irreducible. Thus,
Theorem 2.10 implies

\proclaim{Proposition 2.11}
The tail boundary of an irreducible RWTDI with locally trivial tail boundary
is totally trivial.
\endproclaim

\remark{Remark 2.12}
An immediate generalization of irreducible RWTDI is provided by
{\it periodic\/} RWTDI. We shall say that $\RWTDI(\mmu)$ has period $d\ge 1$ if
there exists a homomorphism $\f:G\to\Z_d$ such that all measures $\mu_n$ are
concentrated on $\f^{-1}(1)$, and the RWTDI on $G_0=\ker\f$ determined by the
sequence of measures
$$
\mu'_n = \mu_{nd}\mu_{nd+1}\cdots\mu_{(n+1)d-1}
$$
is irreducible. In this case $G(\mmu)=G_0$, and the total tail boundary is
isomorphic to $\Z_d$.
\endremark

\subhead
2.4. The entropy
\endsubhead

Below we shall need several facts from the entropy theory of measurable
partitions of Lebesgue spaces \cite{Ro67}. First recall that the
{\it entropy of a discrete probability distribution\/} $p=(p_1,p_2,\dots)$ is
defined as $H(p) = - \sum_i p_i \log p_i$. The {\it entropy $H(\xi)=H_m(\xi)$
of a countable partition\/} $\xi=\{X_i\}$ of a Lebesgue probability space
$(X,m)$ is defined as the entropy of the probability distribution
$p_i=m(X_i)$. In other words,
$$
H(\xi) = - \int \log m(\xi(x))\,dm(x) \;,
$$
where $\xi(x)$ is the element of the partition $\xi$ containing $x$. Given
another measurable (not necessarily countable!) partition $\zeta$ of $X$, the
{\it conditional entropy\/} of $\xi$ with respect to $\zeta$ is defined as
$$
H(\xi|\zeta)
= \int H_{\zeta(x)}(\xi)\,dm(x)
= - \int \log m^{\zeta(x)}( \xi(x)) \,dm(x) \;.
$$
where $x\mapsto\zeta(x)\subset X$ is the projection from the space $(X,m)$ onto
its quotient by the partition $\zeta$ (we identify the points of the quotient
space with the corresponding elements of the partition $\zeta$), the measures
$m^{\zeta(x)}$ are the conditional measures of this projection, and
$H_{\zeta(x)}(\xi)$ is the entropy of $\xi$ with respect to the measure
$m^{\zeta(x)}$.

\proclaim{Proposition 2.13 {\rm (\cite{Ro67}) }}
Let $\xi$ and $\zeta$ be measurable partitions of a Lebesgue space $(X,m)$. If
$\xi$ is countable with $H(\xi)<\infty$, then \roster
\item"(i)"
$0 \le H(\xi|\zeta) \le H(\xi)$, and $H(\xi|\zeta)=0$ (resp.,
$H(\xi|\zeta)=H(\xi)$) iff $\zeta$ is a refinement of $\xi$ (resp., $\zeta$ and
$\xi$ are independent).
\item"(ii)"
If $\zeta'$ is a refinement of $\zeta$, then $H(\xi|\zeta')\le H(\xi|\zeta)$,
and the equality holds iff $m^{\zeta'(x)}(\xi(x))=m^{\zeta(x)}(\xi(x))$ for
$m$-a.e. $x\in X$.
\item"(iii)"
If $\zeta$ is the limit of a monotonously decreasing sequence of measurable
partitions $\zeta_n$, then $H(\xi|\zeta_n)\nearrow H(\xi|\zeta)$.
\endroster
\endproclaim

\proclaim{Theorem 2.14}
If a sequence $\mmu=(\mu_0,\mu_1,\dots)\in \PP G^{\Z_+}$ is such that
$H(\mu_n)<\infty$ for all measures $\mu_n$, then for any $k\in\Z_+$ there
exists a limit
$$
h_k = h_k(\mmu)
= \lim_{n\to\infty} \bigl[ H(\mu_{0,n}) - H(\mu_{k,n}) \bigr] \ge 0 \;,
$$
and the tail boundary of $\RWTDI(\mmu)$ is (locally) trivial iff
$h_k(\mmu)=0$ for all $k$.
\endproclaim

\demo{Proof}
For the coordinate partitions of the path space $(\P,G^{\Z_+})$ the Markov
property yields
$$
H(\a_{0,k}) = \sum_{i=0}^{k-1} H(\mu_i) \;,
\tag 2.4
$$
and for $k<n$
$$
\aligned
H(\a_{0,k}|\a_{n,\infty})
&= \sum_{i=0}^{k-1} H(\mu_i) + H(\mu_{k,n-1}) - H(\mu_{0,n-1}) \\
&= H(\a_{0,k}) + H(\mu_{k,n-1}) - H(\mu_{0,n-1}) \;.
\endaligned
$$
The partitions $\a_{n,\infty}$ are decreasing to the tail
partition $\a_\infty$. Therefore by Proposition 2.13 (i), (iii) the limits
$h_k$ exist, and
$$
H(\a_{0,k}|\a_\infty) = H(\a_{0,k}) - h_k \le H(\a_{0,k}) \;.
\tag 2.5
$$
Moreover, $h_k=0$ iff
$$
H(\a_{0.k}) = H(\a_{0,k}|\a_\infty) \;,
$$
i.e., iff the partitions $\a_{0,k}$ and $\a_\infty$ are independent, see
Proposition 2.13 (i). This is obviously the case if $\a_\infty$ is trivial,
i.e., if the tail boundary is locally trivial. Conversely, if all $h_k$
equal 0, then $\a_\infty$ is independent of all coordinate partitions
$\a_{0,k}$, and therefore of the point partition of the path space, so that
$\a_\infty$ must be trivial.
\demoend

\head
3. Random walks with random transitions probabilities
\endhead

Random walks with random transitions probabilities are a specialization of the
notion of random Markov operators discussed in Section 1.5.

\subhead
3.1. Definitions and preliminaries
\endsubhead

\definition{Definition 3.1}
Let $(\O,\l)$ be a probability Lebesgue measure space endowed with an
invertible ergodic measure preserving transformation $T$, and
$$
\mu:\O\to\PP G \;, \qquad \o\mapsto\mu^\o
$$
be a measurable map. Put
$$
\mmu^\o = \bigl( \mu^\o,\mu^{T\o},\mu^{T^2\o}, \dots \bigr) \in \PP G^{\Z_+}\;.
$$
The family of $\RWTDI(\mmu^\o)$ parameterized by the points $\o\in\O$ is called
a {\it random walk on $G$ with random transition probabilities\/} (RWRTP), for
which we shall use the notation $\RWRTP(\O,\l,T,\mu)$. Below we shall usually
write just $\RWTDI(\o)$ instead of $\RWTDI(\mmu^\o)$.
\enddefinition

Simultaneously with the chains $\RWTDI(\o)$ on $G$ parameterized by the points
$\o\in\O$ we shall also consider the ``global chain'' on the space $\O\times G$
determined by the Markov operator
$$
P:L^\infty(\O\times G,\l\otimes m)\toitself \;, \qquad
Pf(\o,g) = \sum \mu^\o(h) f(T\o,gh) \;,
\tag 3.1
$$
with the transition probabilities $\pi_{\o,g}=\d_{T\o}\otimes g\mu^\o$ (cf.
(1.7)). This chain is homogeneous both in time and space (with respect to the
dissipative action of $G$ on $\O\times G$ by left translations), and the
$\s$-finite measure $\l\otimes m$ is easily seen to be $P$-stationary, so that
the operator $P$ is a {\it covering Markov operator\/} in the sense of
\cite{Ka95}. The corresponding quotient chain on the space $\O$ is
deterministic with the transitions $\o\to T\o$. The projection of the ``global
chain'' starting from a fixed point $\o\in\O$ onto $G$ is the RWTDI determined
by the sequence $\mmu^\o$ (cf. Section 1.5).

Let $\P_{\ov\l}$ be the measure in the path space $(\O\times G)^{\Z_+}$ of the
operator $P$ determined by the initial distribution $\ov\l=\l\otimes\d_e$ on
$\O\times G$. Then by formula (1.9) the map $\Pi:\ov\x\mapsto(\o,\x)$ (1.8) is
an isomorphism between the spaces $\bigl((\O\times G)^{\Z_+},\P_{\ov\l}\bigr)$
and $(\O\times G^{\Z_+},\ov\P)$, where
$$
\ov\P = \int \d_\o\otimes\P_\o \,d\l(\o) \;,
\tag 3.2
$$
and $\P_\o$ denotes the probability measure on the path space $G^{\Z_+}$
corresponding to starting $\RWTDI(\o)$ from the identity of the group.

One can also identify the space of paths $\x\in G^{\Z_+}$ starting from $x_0=e$
with the space $\H\cong G^{\Z_+}$ of increments $\h=(h_0,h_1,\dots)$ by the map
$$
\x\mapsto\h \;,\qquad x_n = h_0 h_1 \cdots h_{n-1} \;, \qquad n>0 \;.
$$
[Although the path space and the space of increments in $G$ both coincide with
$G^{\Z_+}$, their meaning for us is quite different, which is why we use a
separate notation for the space of increments.] Therefore, the space
$(\O\times G^{\Z_+},\ov\P)\cong\bigl((\O\times G)^{\Z_+},\P_{\ov\l}\bigr)$ is
isomorphic to the space $(\O\times\H,\Q)$, where
$$
\Q = \int \d_\o \otimes \bigotimes_{i=0}^\infty \mu^{T^i\o} \,d\l(\o) \;.
$$
Below we shall freely switch between the three descriptions of the same
measure space
$$
(\O\times G^{\Z_+},\ov\P)
\cong \bigl((\O\times G)^{\Z_+},\P_{\ov\l}\bigr)
\cong (\O\times\H,\Q)
$$
using the correspondence
$$
\aligned
&(\o,\x) = \bigl(\o,(e,x_1,x_2,\dots) \bigr) \\
&\;\lra\; \ov\x
= (\ov x_n)
= \bigl( (\o,e),(T\o,x_1),(T^2\o,x_2),\dots \bigr) \\
&\;\lra\; (\o,\h)
= \bigl(\o,(h_n)\bigr)
= \bigl(\o,(x_1,x_1^{-1}x_2,x_2^{-1}x_3,\dots)\bigr) \;.
\endaligned
\tag 3.3
$$
For uniformity we shall usually refer to this measure space as
$(\O\times G^{\Z_+},\ov\P)$ ``changing variables'' by formulas (3.3) if
necessary.

Denote by $E$ (resp., $\G$) the tail (resp., the Poisson) boundary of the
operator $P$ (3.1). Below we shall only be interested in the boundary behaviour
of the operator $P$ and of $\RWTDI(\o)$ which is not reducible to a dependence
on the starting point (cf. Section 2.3). Therefore, the tail boundary $E$ will
always be endowed with the tail measure $\e=\tail\ov\P$. By
$\nu=\bnd\ov\P=\p_\G\e$ (where $\p_\G$ is the projection $E\to\G$) we denote
the corresponding harmonic measure on the Poisson boundary.

The fibers of the projection $\p_\O:(E,\e)\to(\O,\l)$ are the tail boundaries
$(E_\o,\e_\o)$ of $\RWTDI(\o)$, where $\e_\o$ is the tail measure on $E_\o$
corresponding to starting $\RWTDI(\o)$ from the group identity at time $0$
(Proposition 1.11). More generally, by $\e_{n,\o}$ we denote the tail measure
on $E_\o$ corresponding to starting $\RWTDI(\o)$ from the group identity at an
arbitrary time $n\ge 0$. In view of Proposition 1.11 we may also consider the
measures $\e_\o,\e_{n,\o}$ as measures on $E$. By
$\nu_\o=\p_\G\e_\o=\p_\G\e_{n,\o}$ we denote the corresponding harmonic
measures on the Poisson boundary $\G$ of the operator $P$. In other words,
$\e_\o=\tail(\P_\o)$ and $\nu_\o=\bnd(\P_\o)$.

The action of $G$ descends from the path space to the tail boundary and to the
Poisson boundary, and clearly $g\e_o=\tail(g\P_\o)$ and $g\nu_\o=\bnd(g\P_\o)$
are the measures on $E$ (resp., on $\G$) corresponding to starting
$\RWTDI(\o)$ from the point $g\in G$ at time $0$. Then formula (2.1) takes in
this setup the form
$$
\e_{n,\o} = \mu^\o \e_{n+1,T\o} \;.
$$
Therefore,
$$
\nu_\o = \mu^\o \nu_{T\o} \;.
$$
The Poisson formulas for bounded harmonic sequences and functions of the
operator $P$ (see Theorems 1.7 and 1.10) take the form, respectively,
$$
f_n (\o,g) = \la \wh f, g\e_{n,\o} \ra
$$
and
$$
f (\o,g) = \la \wh f, g\nu_\o \ra \;.
$$

\proclaim{Proposition 3.2}
The tail boundaries of the $\RWTDI(\o),\;\o\in\O$ are all trivial or
non-trivial simultaneously.
\endproclaim

\demo{Proof}
As it follows from Theorem 2.6 (see Remark 2.7), triviality of the
tail boundary of $\RWTDI(\o)$ implies triviality of the tail boundary of
$\RWTDI(T\o)$. Therefore, the claim follows from ergodicity of the
transformation $T$.
\demoend

We shall say that $\RWRTP(\O,\l,T,\mu)$ is {\it irreducible\/} if $\RWTDI(\o)$
is irreducible for $\l$-a.e. $\o\in\O$ (see Definition 2.9). Proposition 1.11
then implies

\proclaim{Proposition 3.3}
If $\RWRTP(\O,\l,T,\mu)$ is irreducible and the tail boundary of $\l$-a.e.
$\RWTDI(\o),\;\o\in\O$  is locally trivial, then the Poisson boundary
$\G$ of the operator $P$ is trivial.
\endproclaim

\remark{Remark 3.4}
Proposition 3.3 remains true if $\RWRTP(\O,\l,T,\mu)$ is {\it periodic\/},
i.e., if a.e. $\RWTDI(\o)$ is periodic with the same period $d$.
\endremark

\proclaim{Problem 3.5}
Give general conditions which would ensure stability (in the sense of
Definition 1.12) of $\RWRTP(\O,\l,T,\mu)$.
\endproclaim

\subhead
3.2. The asymptotic entropy
\endsubhead

\definition{Definition 3.6}
We shall say that $\RWRTP(\O,\l,T,\mu)$ has {\it finite entropy\/} if
$$
\ov H = \ov H(\O,\l,T,\mu) = \int H(\mu^\o)\,d\l(\o) < \infty \;.
$$
\enddefinition

\proclaim{Theorem 3.7}
If $\RWRTP(\O,\l,T,\mu)$ has finite entropy, then the limit
$$
\ov h = \ov h(\O,\l,T,\mu) = \lim_{n\to\infty} \frac{H(\mu_{0,n-1}^\o)}n
\tag 3.4
$$
exists a.e. and in the space $L^1(\O,\l)$ and is independent of $\o$.
\endproclaim

\definition{Definition 3.8}
The limit (3.4) is called the {\it asymptotic entropy\/} of
$\RWRTP(\O,\l,T,\mu)$.
\enddefinition

\demo{Proof of Theorem 3.7}
The measure $\mu_{0,k+n-1}^\o$ is the convolution of the measures
$\mu_{0,k-1}^\o$ and $\mu_{k,k+n-1}^\o=\mu_{0,n-1}^{T^k\o}$. Therefore,
$$
H(\mu_{0,k+n-1}^\o) \le H(\mu_{0,k-1}^\o) + H(\mu_{0,n-1}^{T^k\o})
$$
so that the sequence of functions $\f_n(\o)=H(\mu_{0,n-1}^\o)$ on $\O$
satisfies conditions of the Kingman subadditive ergodic theorem
(e.g., see \cite{De80}), which implies the claim.
\demoend

\proclaim{Theorem 3.9}
If $\ov H(\O,\l,T,\mu)<\infty$, then $\ov h(\O,\l,T,\mu)=0$ iff the tail
boundary of $\l$-a.e. $\RWTDI(\o)$ is trivial.
\endproclaim

\proclaim{Corollary}
If $\ov h(\O,\l,T,\mu)=0$, then the Poisson boundary of $\RWRTP(\O,\l,T,\mu)$
is trivial.
\endproclaim

\remark{Remark 3.10}
Contrary to the situation with the ``ordinary'' random walks we do not know
whether triviality of the Poisson boundary of $\RWRTP(\O,\l,T,\mu)$ implies
that\linebreak
 $\ov h(\O,\l,T,\mu)=0$. The reason for this difference is that for
ordinary random walks the tail and the Poisson boundary coincide with respect
to any single point initial distribution (e.g., see \cite{Ka92}), so that the
entropy criterion of triviality of the tail boundary automatically becomes the
entropy criterion of triviality of the Poisson boundary. However, for RWRTP the
relation between the tail and the Poisson boundaries is more complicated, see
Definition 1.12 (together with the discussion there) and Problem 3.5. Of
course, if the Poisson boundary is stable in the sense of Definition 1.12, then
it is trivial iff the tail boundary $\G$ coincides with $\O$, i.e., iff $\ov
h(\O,\l,T,\mu)=0$.
\endremark

\demo{Proof of Theorem 3.9}
By Theorem 2.14 for $\l$-a.e. $\o\in\O$ and any $k>0$ there exists the limit
$$
h(\o)
= \lim_{n\to\infty} \bigl[ H(\mu_{0,n}^\o) - H(\mu_{1,n}^\o) \bigr]
= \lim_{n\to\infty} \bigl[ H(\mu_{0,n}^\o) - H(\mu_{0,n-1}^{T\o}) \bigr] \;,
\tag 3.5
$$
and the tail boundary of $\l$-a.e. $\RWTDI(\o)$ is locally trivial iff
$h(\o)=0$ for $\l$-a.e. $\o\in\O$. Since $h(\o)\le H(\mu^\o)$, the function
$h$ is integrable, and the convergence in (3.5) also holds in the space
$L^1(\O,\o)$, whence
$$
\int h(\o)\,d\l(\o)
= \lim_{n\to\infty}
\left[ \int H(\mu_{0,n}^\o)\,d\l(\o) - \int H(\mu_{0,n-1}^\o)\,d\l(\o) \right]
\;.
$$
Therefore, by (3.4)
$$
\int h(\o)\,d\l(\o) = \ov h(\O,\l,T,\mu) \;.
\tag 3.6
$$
In particular, $\ov h(\O,\l,T,\mu)=0$ iff $h\equiv 0$.
\demoend

\definition{Definition 3.11 {\rm (\cite{Ka98})}}
A probability measure $\Th$ on $G^{\Z_+}$ has {\it asymptotic
entropy\/} $\hh(\Th)$ if the following Shannon--Brei\-man--McMillan type
{\it equidistribution condition\/} is satisfied:
$$
- \frac1n \log \th_n(x_n)\to \hh(\Th)
$$
for $\Th$-a.e. sequence $\x=\{x_n\}\in G^{\Z_+}$ and in the space
$L^1(G^{\Z_+},\Th)$, where $\th_n$ are the one-dimensional distributions of the
measure $\Th$.
\enddefinition

The following result shows that for RWRTP the asymptotic entropies in the sense
of Definitions 3.8 and 3.11 coincide.

\proclaim{Theorem 3.12}
If $\ov H(\O,\l,T,\mu)<\infty$, then
$$
\hh(\P_\o) = \ov h(\O,\l,T,\mu)
$$
for $\l$-a.e. $\o\in\O$.
\endproclaim

In combination with Theorem 3.9 it immediately implies

\proclaim{Corollary}
If $\ov H(\O,\l,T,\mu)<\infty$ and for $\l$-a.e. $\o\in\O$ there exists a
sequence of finite subsets $A_n\subset G$ such that $\log|A_n|=o(n)$ and
$\mu^\o_{0,n}(A_n)>\e$ for a certain fixed number $\e>0$, then the tail
boundary of $\RWTDI(\o)$ is trivial for $\l$-a.e. $\o\in\O$.
\endproclaim

\proclaim{Lemma 3.13}
Let $(S\h)_n=h_{n+1}$ be the shift in the space of increments $\H$. Then
the transformation
$$
\ov T (\o,\h) = (T\o,S\h)
$$
of the space $(\O\times G^{\Z_+},\ov\P)\cong (\O\times\H,\Q)$ (see {\rm (3.3)})
is measure preserving and ergodic.
\endproclaim

\demo{Proof}
The map
$$
\bigl(\o,(h_0,h_1,\dots)\bigr) \mapsto \bigl((\o,h_0),(T\o,h_1),\dots\bigr)
$$
identifies the space $(\O\times\H,\Q)$ with the path space of the Markov chain
on $\O\times G$ with the transition probabilities
$$
\Prob \bigl[ (T\o,h') | (\o,h) \bigr] = \mu^{T\o}(h')
\tag 3.7
$$
and the initial distribution $d\th(\o,h)=d\l(\o)\mu^\o(h)$. This identification
conjugates $\ov T$ with the shift in $(\O\times G)^{\Z_+}$. By $T$-invariance
of $\l$ the measure $\th$ is stationary with respect to the transition
probabilities (3.7), so that the measure $\Q$ is $\ov T$-invariant (this fact
can be also easily checked directly).

Further, since the measure $\th$ is finite, ergodicity of $\Q$ with respect to
$\ov T$ is equivalent to absence of non-constant bounded harmonic functions of
the chain (3.7), e.g., see \cite{Ka92}. The transition probabilities (3.7) from
a point $(\o,h)$ do not depend on $h$, so that any such function depends on
$\o$ only, i.e., is a non-constant $T$-invariant function on $\O$, which is
impossible by ergodicity of $T$.
\demoend

\remark{Remark 3.14}
In terms of the path space $(\O\times G)^{\Z_+}$ the transformation $\ov T$
takes the form
$$
\bigl( (\o,e),(T\o,x_1),(T^2\o,x_2), \dots \bigr)
\mapsto \bigl( (T\o,e),(T^2\o,x_1^{-1}x_2),(T^3\o,x_1^{-1}x_3),\dots \bigr) \;,
$$
i.e., it is the combination of the shift in the path space with the subsequent
group translation consisting in moving the origin of the shifted path to the
identity of $G$.
\endremark

\demo{Proof of Theorem 3.12}
Put
$$
\f_n(\o,\h) = - \log \mu_{0,n-1}^\o(x_n)
= - \log \mu_{0,n-1}^\o(h_0 h_1 \cdots h_{n-1}) \;.
$$
Since
$$
\mu^\o_{0,k+n-1}(h_0 h_1 \cdots h_{k+n-1})
\ge \mu^\o_{0,k-1}(h_0 h_1 \cdots h_{k-1})
\mu^\o_{k,k+n-1}(h_k h_{k+1} \cdots h_{k+n-1}) \;,
$$
we have
$$
\f_{k+n}(\o,\h) \le \f_k(\o,\h) + \f_n(T^k\o,S^k\h) \;.
$$
Finiteness of entropy of $\RWRTP(\O,\l,T,\mu)$ means that
$$
\aligned
\int \f_1(\o,\h)\,d\Q(\o,\h)
&= - \int \log \mu^\o(h_0)\,d\Q(\o,\h) \\
&= - \int \sum_g \log\mu^\o(g)\,\mu^\o(g)\,d\l(\o) \\
&= \int H(\mu^\o)\,d\l(\o) = \ov H(\O,\l,T,\mu) < \infty \;,
\endaligned
\tag 3.8
$$
so that the conditions of the Kingman subadditive ergodic theorem are
satisfied, and therefore there exists a limit
$$
\lim_{n\to\infty} -\frac1n \log\mu_{0,n-1}^\o(h_0h_1\cdots h_{n-1})
\tag 3.9
$$
for $\Q$-a.e. pair $(\o,\h)\in\O\times\H$ and in the space
$L^1(\O\times\H,\Q)$. Since
$$
\int \f_n(\o,\h)\,d\Q(\o,\h)
= \int H(\mu^\o_{0,n-1})\,d\l(\o)
$$
(cf. (3.8)), and
$$
\frac1n \int H(\mu^\o_{0,n-1})\,d\l(\o) \to \ov h(\O,\l,T,\mu)
$$
by Theorem 3.7, the limit (3.9) coincides with $\ov h(\O,\l,T,\mu)$.
\demoend

\subhead
3.3. The asymptotic entropy of conditional chains
\endsubhead

We shall apply to the partitions of the path space of the operator $P$ (3.1)
the notations $\ov\a_k, \ov\a_{k,l}$, etc. introduced in Section 1.2 (sometimes
we overline the objects connected with the operator $P$ on the state space
$\O\times G$ in order to distinguish them from the objects associated with the
individual $\RWRTP(\o)$ on $G$). We continue to use the identifications (3.3).
In the model $(\O\times\H,\Q)$ of the path space $((\O\times
G)^{\Z_+},\P_{\ov\l})\cong(\O\times G^{\Z_+},\ov\P)$ the partition
$\ov\a_{0,n}$ coincides with the common refinement of the point partition of
the space $\O$ and the partition of $\H$ determined by the first $n$
coordinates $h_0,h_1,\dots,h_{n-1}$. In particular, $\ov\a_0$ coincides
$\ov\P$-mod 0 with the preimage partition of the projection
$\Pi_\O:\ov\x\mapsto\o$ (1.10). By $\ov\a_\infty=\lim\ov\a_{n,\infty}$ we
denote the tail partition. Note that $\ov\a_\infty$ is a refinement of
$\ov\a_0$, see Proposition 1.11.

Below all the entropies and the conditional entropies of partitions of the
space $\O\times G^{\Z_+}$ are calculated with respect to the measure $\ov\P$.
We begin with expressing the asymptotic entropy $\ov h(\O,\l,T,\mu)$ in terms
of the conditional entropies of coordinate partitions. Integrating formula
(2.5) with respect to the measure $\l$ and using (2.4), (3.5) and (3.6), we
obtain

\proclaim{Lemma 3.15}
If $\ov H(\O,\l,T,\mu)<\infty$, then for any $k\ge 0$
$$
H(\ov\a_{0,k}|\ov\a_\infty) = k \Bigl[ \ov H(\O,\l,T,\mu) - \ov h(\O,\l,T,\mu)
\Bigr] \;.
$$
\endproclaim

By $\s_E$ (resp., $\s_\O$) we denote the point partition of the tail boundary
$(E,\e)$ (resp., the partition generated by the projection $\p_\O:E\to\O$, see
Proposition 1.11). Let $\xi$ be a $G$-invariant partition of $E$ which is a
refinement of $\s_\O$. Denote by $\ov\a_\xi=\tail^{-1}\xi$ the partition of the
path space $((\O\times G)^{\Z_+},\P_{\ov\l})\cong(\O\times G^{\Z_+},\ov\P)$
which is the preimage of $\xi$ under the map $\tail$. Since $\tail^{-1}\s_E$ is
the tail partition $\ov\a_\infty$, and $\tail^{-1}\s_\O=\ov\a_0$, we have
$$
\ov\a_0 \precc \ov\a_\xi \precc \ov\a_\infty \;.
$$
The quotient $(E^\xi,\e^\xi)$ of the tail boundary $(E,\e)$ by the partition
$\xi$ can be also considered as the quotient of the path space by the partition
$\ov\a_\xi$. Denote by
$$
\tail^\xi:(\O\times G^{\Z_+},\ov\P)\to (E^\xi,\e^\xi)
$$
the associated quotient map. The spaces $(E^\xi,\e^\xi)$ are ``random
analogues'' of the $\mu$-boundaries in the case of usual time homogeneous
random walks on groups, see \cite{Fu71}, \cite{Ka00}. If $\xi=\s_E$, then
$\tail^\xi=\tail$, and if $\xi=\s_\O$, then $\tail^\xi$ coincides with the
projection $\Pi_\O$ (1.10). Denote by $\e_\o^\xi=\e_{0,\o}^\xi$ (resp.,
$\e_{n,\o}^\xi$) the images of the tail measures $\e_\o$ (resp., $\e_{n,\o}$)
under the projection $E\to E^\xi$.

\proclaim{Lemma 3.16}
If $\ov H(\O,\l,T,\mu)<\infty$, then for any $k\ge 0$
$$
\align
H(\ov\a_{0,k}|\ov\a_\xi)
&= k H(\ov\a_{0,1}|\ov\a_\xi) \\
&= k \left[\ov H(\O,\l,T,\mu)
- \int \log\frac {dx_1\e_{1,T\o}^\xi} {d\e_\o^\xi}
(\tail^\xi(\o,\x))\,d\ov\P(\o,\x)\right] \;.
\endalign
$$
\endproclaim

\demo{Proof}
Since $\s_\O \precc \xi$, conditioning by $\xi$ uniquely determines the
starting point $\o$ of the sample path $\ov\x\lra (\o,\x)$. The traces of $\xi$
on the elements of the partition $\s_\O$ (i.e., on the tail boundaries $E_\o$,
see Proposition 1.11) are $G$-invariant partitions, so that we may apply
Proposition 2.4, according to which the conditional probability of the element
of the partition $\ov\a_{0,k}$ containing given $(\o,\x)\lra (\o,\h)$ with
respect to the partition $\ov\a_\xi$ is
$$
\mu^\o(h_1)\mu^{T\o}(h_2)\dots \mu^{T^{k-1}\o}(h_k)
\frac{dx_k\e_{k,T^k\o}^\xi}{d\e_\o^\xi}(\tail^\xi(\o,\x)) \;,
$$
whence integrating we obtain
$$
H(\ov\a_{0,k}|\ov\a_\xi)
= k \ov H(\O,\l,T,\mu)
- \int\log \frac {dx_k \e_{k,T^k\o}^\xi}
{d\e_\o^\xi}(\tail^\xi(\o,\x))\,d\P(\o,\x) \;.
$$
The integrand in the last term in the right hand side telescopes as
$$
\log \frac{dx_k\e_{k,T^k\o}^\xi}{d\e_\o^\xi}(\tail^\xi(\o,\x))
= \sum_{i=0}^{k-1} \f(\ov T^i(\o,\x)) \;,
$$
where
$$
\f(\o,x) = \log \frac{dx_1\e_{1,T\o}^\xi}{d\e_\o^\xi}(\tail^\xi(\o,\x)) \;,
$$
and $\ov T$ is the transformation of the path space $\O\times G^{\Z_+}$
introduced in Lemma 3.13. Since $\ov T$ preserves the measure $\ov\P$, we get
the claim.
\demoend

\proclaim{Corollary}
If $\ov H(\O,\l,T,\mu)<\infty$, then
$$
\ov h(\O,\l,T,\mu) = \int \log\frac {dx_1\e_{1,T\o}} {d\e_\o}
(\tail(\o,\x)) \;.
$$
\endproclaim

\proclaim{Lemma 3.17}
Let $\xi$ and $\xi'$ be two $G$-invariant measurable partitions of the tail
boundary $E$ such that $\s_\O\precc\xi\precc\xi'$. If
$\ov H(\O,\l,T,\mu)<\infty$, then
$$
H(\ov\a_{0,1}|\ov\a_\xi)\ge H(\ov\a_{0,1}|\ov\a_{\xi'}) \;,
$$
and the equality holds iff $\xi=\xi'$.
\endproclaim

\demo{Proof}
Obviously, if $\xi'$ is a refinement of $\xi$, then $\ov\a_{\xi'}$ is a
refinement of $\ov\a_\xi$, so that the inequality follows from Proposition 2.13
(ii). If $H(\ov\a_{0,1}|\ov\a_\xi)=H(\ov\a_{0,1}|\ov\a_{\xi'})$, then by Lemma
3.16, $H(\ov\a_{0,k}|\ov\a_\xi)=H(\ov\a_{0.k}|\ov\a_{\xi'})$ for any $k\ge 1$,
which by Proposition 2.13 (ii) means that for $\e$-a.e. point $\g\in E$ the
$k$-dimensional distributions of the conditional measures $\ov\P^{\xi(\g)}$ and
$\ov\P^{\xi'(\g)}$ are the same. Therefore, for $\e$-a.e. $\g\in E$ the
conditional measures $\ov\P^{\xi(\g)}$ and $\ov\P^{\xi'(\g)}$ coincide, which
is only possible if $\xi=\xi'$.
\demoend

\proclaim{Corollary}
Let $\xi$ be a $G$-invariant measurable partitions of the tail boundary $E$
such that $\s_\O\precc\xi$. Then $\xi=\s_E$ iff
$$
H(\a_{0,1}|\a_\xi) = H(\a_{0,1}|\a_\infty) \;,
$$
\endproclaim

\proclaim{Theorem 3.18}
Let $\xi$ be a $G$-invariant measurable partition of the tail boundary
$E$ such that $\s_\O\precc\xi$. If $\ov H(\O,\l,T,\mu)<\infty$, then for
$\e^\xi$-a.e. point $\xi(\g)\in E^\xi$ the asymptotic entropy (in the sense of
Definition 3.11) of the conditional measure $\ov\P^{\xi(\g)}$ exists and is
equal to
$$
\hh(\ov\P^{\xi(\g)})
= \ov h(\O,\l,T,\mu) - \int \log\frac {dx_1\e_{1,T\o}^\xi} {d\e_\o^\xi}
(\tail^\xi(\o,\x))\,d\ov\P(\o,\x) \;.
\tag 3.10
$$
\endproclaim

\demo{Proof}
Since $\s_\O\precc\xi$, by Proposition 2.4 the one-dimensional distributions
$\pi_n^{\xi(\g)}$ of the conditional measure $\ov\P^{\xi(\g)}$ are
$$
\pi_n^{\xi(\g)}(g)
= \mu^\o_{0,n-1}(g) \frac{dg\e^\xi_{n,T^n\o}}{d\e^\xi_\o}(\xi(\g)) \;,
$$
where the point $\o=\o(\g)\in\O$ is determined by the projection $E^\xi\to\O$.
Theorem 3.12 implies the convergence ($\ov\P$-a.e. and in the space
$L^1(\O\times G^{\Z_+},\ov\P)$)
$$
-\frac1n \log \mu^\o_{0,n-1}(x_n) \to \ov h(\O,\l,T,\mu) \;,
\tag 3.11
$$
and the telescoping at the end of the proof of Lemma 3.16 in combination
with the Birkhoff ergodic theorem for the transformation $\ov T$ yields the
convergence (also a.e. and in the $L^1$-space),
$$
\frac1n \log \frac{dx_n\e^\xi_{n,T^n\o}}{d\e^\xi_\o}(\tail^\xi(\o,\x))
\to \int \log\frac {dx_1\e_{1,T\o}^\xi} {d\e_\o^\xi}
(\tail^\xi(\o,\x))\,d\ov\P(\o,\x) \;.
\tag 3.12
$$

Combining (3.11) and (3.12), we obtain the convergence ($\ov\P$-a.e. and in the
 $L^1$-space) of
$$
-\frac1n\log\pi_n^{\tail^\xi(\o,\x)}(x_n)
$$
to the limit (3.10), which implies the claim, because the measures
$\ov\P^{\xi(\g)}$ are the conditional measures of the measure $\ov\P$ with
respect to the partition $\ov\a_\xi$.
\demoend

Now, combining Theorem 3.18 with Lemmas 3.15, 3.16 and 3.17, we get the
following generalization of Theorem 3.9

\proclaim{Theorem 3.19}
Let $\xi$ be a $G$-invariant measurable partition of the tail boundary
$E$ such that $\s_\O\precc\xi$. If $\ov H(\O,\l,T,\mu)<\infty$, then $\xi=\s_E$
iff the asymptotic entropies of the conditional measures $\ov\P^{\xi(\g)}$
vanish.
\endproclaim

\proclaim{Corollary}
The partition $\xi$ coincides with $\s_E$ iff for $\e^\xi$-a.e. point
$\xi(\g)\in E^\xi$ there exist $\e>0$ and a sequence of sets
$A_n=A_n(\xi(\g))\subset G$ such that $\log|A_n|=o(n)$ and
$\pi^{\xi(\g)}_n(A_n)>\e$ for all sufficiently large $n$.
\endproclaim

\head
4. Triviality and identification of the tail and the Poisson boundaries
\endhead

In this section we consider several concrete classes of groups and describe the
boundaries of RWRTP on these groups.

\subhead
4.1. Boundary triviality
\endsubhead

The entropy theory developed in Section 3 allows one to extend to RWRTP almost
all results on triviality and identification of the boundaries earlier
obtained for usual random walks, see \cite{KV83}, \cite{Ka00}.

Throughout this section we assume that the group $G$ acts by isometries
on a complete metric space $(X,d)$. Fix once and forever a reference point
$o\in X$ (its choice is irrelevant for what follows) and put
$$
|g| = |g|_X = d(o,go) \;,\qquad g\in G \;.
$$
Suppose that the group $G$ has {\it bounded exponential
growth\/} with respect to the space $X$, i.e.,
$$
v(G,X) = \limsup_{t\to\infty} \frac1t \log\card\{g\in G: |g|_X \le t\}
< \infty \;.
\tag 4.1
$$
If $G$ is a finitely generated group, and $X\cong G$ is its Cayley graph
determined by a finite generating set and endowed with the associated word
metric $d$, then condition (4.1) is obviously satisfied. Another example is
provided by a discrete subgroup $G$ of isometries of a Riemannian manifold of
bounded geometry $X$. If $v(G,X)=0$ we shall say that the group $G$ has
{\it subexponential growth\/}.

For a measure $\th\in\PP G$ denote by
$$
|\th| = \sum_g d(o,go) \th(g)
$$
its {\it first moment\/}. We shall say that $\RWRTP(\O,\l,T\,\mu)$ has a
{\it finite first moment\/} (with respect to the space $X$) if
$$
\int |\mu^\o|\,d\l(\o) < \infty \;.
$$

Using the triangle inequality and the Kingman subadditive ergodic theorem, we
derive

\proclaim{Theorem 4.1}
If $\RWRTP(\O,\l,T,\mu)$ on the group $G$ has a finite first moment with
respect to the space $X$ then there exists a number $l=l(\O,\l,T,\mu,X)$
called the {\it linear rate of escape\/} such that for $\ov\P$-a.e.
$(\o,\x)\in\O\times G^{\Z_+}$,
$$
\lim_{n\to\infty} \frac{|x_n|_X}n = l \;.
$$
The convergence also holds in the space $L^1(\O\times G^{\Z_+},\ov\P)$, where
$\ov\P$ is the measure {\rm (3.2)}.
\endproclaim

\proclaim{Lemma 4.2 {\rm (\cite{De86}) }}
There exists a constant $C=C(G,X)$ such that for any measure $\th\in\PP G$,
$$
H(\th) \le C(|\th|_K + 1) \;.
$$
\endproclaim

Now, using Theorem 3.12 we obtain in the same way as for ordinary
random walks on groups (see \cite{Gu80}) the following result

\proclaim{Theorem 4.3}
If $\RWRTP(\O,\l,T,\mu)$ on $G$ has a finite first moment, then its entropy
$\ov h(\O,\l,T,\mu)$, the rate of escape $l(\O,\l,T,\mu,X)$ and the rate of
growth $v(G,X)$ satisfy the inequality
$$
\ov h \le l v \;.
$$
\endproclaim

\proclaim{Corollary {\rm (cf. Proposition 2.9)} }
If $\RWRTP(\O,\l,T,\mu)$ is irreducible and the group $G$ is non-amenable,
then $l>0$.
\endproclaim

Theorem 4.3 in combination with Theorem 3.9 implies triviality of the tail
boundaries of $\l$-a.e. $\RWTDI(\o)$ and of the Poisson boundary of
$\RWRTP(\O,\l,T,\mu)$ when either $l(\O,\l,T,\mu,X)$ or $v(G,X)$ vanish. Since
any finitely generated nilpotent group has polynomial growth (with respect to
the word metric determined by any finite generating set), we obtain

\proclaim{Theorem 4.4}
The Poisson boundary of any RWRTP with a finite first moment on a finitely
generated nilpotent group is trivial.
\endproclaim

Combining Theorem 4.4 with Proposition 2.9 we now obtain

\proclaim{Theorem 4.5}
Let $(\O,\l,T,\mu)$ be an irreducible RWRTP with a finite first moment on a
finitely generated nilpotent group $G$. Then for any $g\in G$ and $\l$-a.e.
$\o\in\O$
$$
\| g\mu_{0,n}^\o - \mu_{0,n}^\o \| \toto_{n\to\infty} 0 \;,
$$
i.e., a.e. sequence $\mu_{0,n}^\o$ strongly converges to left invariance on
$G$.
\endproclaim

\remark{Remark 4.6}
By completely different methods Theorem 4.5 was proved in [MR94], [LRW94], and
[Ru95] for compact and abelian groups without any additional moment
assumptions.
\endremark

Returning to Theorem 4.3, recall that another way of proving boundary
triviality consists in showing that the rate of escape $l(\O,\l,T,\mu,X)$
vanishes. The methods used in \cite{Ka91} allow one to do it for ``centered''
RWRTP on several classes of solvable groups in the same way as for usual time
homogeneous random walks. For the sake of brevity we shall consider just the
class of {\it polycyclic groups\/}. Without loss of generality we may assume
that the polycyclic group $G=A\sd N$ is the semi-direct product of an abelian
group $A\cong\Z^d$ and a normal finitely generated nilpotent subgroup $N$ (see
\cite{Ka91}). For a measure $\th$ on $G$ denote by $\th_A$ its projection onto
$A$, and by
$$
\ov\th_A = \sum_{a\in A} \th_A(a) a \in \R^d
$$
the barycenter of $\th_A$ (this definition requires finiteness of the first
moment of the measure $\th_A$). If $\ov\th_A=0$, then the measure $\th$ is
called {\it centered\/}. We shall say that $\RWRTP(\O,\l,T,\mu)$ on $G$ with a
finite first moment is centered if
$$
\int\ov\mu^\o_A\,d\l(\o) < \infty \;.
$$

\proclaim{Theorem 4.7}
The Poisson boundary of any centered RWRTP with a finite first moment on a
polycyclic group is trivial.
\endproclaim

\proclaim{Theorem 4.8} For any centered irreducible RWRTP with a finite first
moment on a polycyclic group $G$ a.e. sequence $\mu_{0,n}^\o$ strongly
converges to left invariance on $G$.
\endproclaim

\subhead
4.2. Boundary identification
\endsubhead

We shall now look at the problem of identifying the tail and the Poisson
boundaries of RWRTP on groups. Suppose, for the sake of argument, that our
group $G$ admits an invariant compactification $\ov G$ with the boundary $\part
G$ (i.e., the action of $G$ on itself by left translations extends to a
continuous action on $\ov G$), and that $\P_\o$-a.e. sample path $\x=(x_n)$
converges in this compactification to a limit point
$x_\infty=x_\infty(\x)\in\part G$ for $\l$-a.e. $\o\in\O$. Obviously, the map
$\x\mapsto x_\infty$ is measurable with respect to the tail $\s$-algebra of the
global operator $P$ on $\O\times G$ (actually, the topological nature of $\part
G$ is completely irrelevant for what follows). Therefore, the map
$(\o,\x)\mapsto (\o,x_\infty)$ gives rise to a measurable partition $\xi$ of
the Poisson boundary $E$ of the operator $P$. The partition $\xi$ is
$G$-invariant, and it is a refinement of the partition $\s_\O$ determined by
the projection $E\to\O$. Coincidence of the partition $\xi$ with the point
partition $\s_E$ of $E$ means that the tail boundary $E$ actually can be
identified with the product $\O\times\part G$. Therefore, in the latter case
the Poisson boundary of the RWRTP is stable (in the sense of Definition~1.12),
i.e., the Poisson boundary of the operator $P$ and the tail boundaries $E^\o$
can be both identified with $\part G$.

The main method of proving boundary convergence for ``groups with hyperbolic
properties'' goes back to Furstenberg \cite{Fu71} and consists in using the
martingale convergence theorem in combination with contracting (proximality)
properties of the action of $G$ on the boundary $\part G$, see \cite{CS89},
\cite{Wo93}, \cite{Ka00}. This method does not impose any moment conditions on
the random walk and may be combined with the ``strip approximation'' criteria
\cite{Ka00} to give a full boundary identification. However, its application
to RWRTP is rather tedious and we could not get rid of rather awkward
conditions on the measures $\mu^\o$ (like existence of a single
non-degenerate measure on $G$ dominated by a.e. $\mu^\o$) following this way.
Instead of this we shall use the ``ray approximation'' approach (see
\cite{Ka00}) and its recent generalization obtained in
\cite{KMa99} which will save us from a good deal of technical details.

Recall that a metric space $(X,d)$ is called convex if for any two points
$x,y\in X$ there exists a {\it midpoint\/} $z\in X$ such that
$$
d(x,z) = d(y,z) = \frac12 d(x,y) \;.
$$
In a complete convex space any two points can be joined by a geodesic (see
the related definitions in [BH99]).

A metric space $(X,d)$ is called {\it uniformly convex\/} if it is convex and
in addition there exists a strictly decreasing continuous function $\f$ on
$[0,1]$ with $\f(0)=1$ such that for any $x,y,w\in X$ and a midpoint
$z$ of $x$ and $y$
$$
\frac{d(z,w)}R \le \f\left(\frac{d(x,y)}{2R}\right) \;,
$$
where $R=\max\{d(x,w),d(y,w)\}$. The midpoints (and therefore geodesics with
given endpoints) in a uniformly convex space are unique.

A convex metric space $(X,d)$ is called {\it non-positively curved\/} (in the
sense of Busemann) if for any $x,y,z\in X$ and any midpoints $m_{xz}$ (resp.,
$m_{yz}$) of $x$ and $y$ (resp., of $y$ and $z$)
$$
d(m_{xz}.m_{yz})\le\frac12 d(x,y) \;.
$$

From now on we shall assume that
$$
\ft{The metric space $X$ on which the group
$G$ acts  is uniformly convex and satisfies the Busemann non-positive curvature
condition.}
\tag 4.2
$$
Denote by $\part X$ the space of asymptotic classes of geodesic rays in $X$. We
shall identify $\part X$ with the space of geodesic rays issued from the point
$o$. Examples of spaces $(X,d)$ satisfying condition (4.2) include all
Cartan--Hadamard manifolds (in particular, non-compact Riemannian symmetric
spaces without compact factors) and all metric trees. In the first case
$\part X$ is the visibility sphere of $X$, and in the second case it is the
space of ends of $X$.

An application of \cite{KMa99} to the transformation $\ov T$ of the space
$(\O\times G^{\Z_+},\ov\P)$ gives

\proclaim{Theorem 4.9}
Suppose that $\RWRTP(\O,\l,T,\mu)$ on the group $G$ has a finite first moment,
its rate of escape $l=l(\O,\l,T,\mu,X)$ is positive, and the space $X$
satisfies conditions (4.1) and (4.2). Then for
$\ov\P$-a.e. $(\o,\x)\in\O\times G^{\Z_+}$ there exists a unique geodesic ray
$\g=\g(\o,\x)\in\part X$ such that
$$
d(x_n,\g(nl))=o(n) \;.
$$
\endproclaim

Note that in view of the Corollary of Theorem 4.3 the condition $l>0$ in
Theorem 4.9 is not really restrictive as discrete groups of isometries of
non-positively curved spaces are usually non-amenable. Theorem 4.9 in
combination with Theorem 3.19 immediately implies

\proclaim{Theorem 4.10}
Under conditions of Theorem 4.9 the Poisson boundary of RWRTP is stable and is
isomorphic to the space $\part X$ with the resulting hitting measure.
\endproclaim

Therefore, the Poisson boundary identifies with the natural geometric
boundaries for RWRTP with a finite first moment on free groups and on discrete
groups of isometries of Cartan--Hadamard manifolds (in particular, in discrete
subgroups of semi-simple Lie groups), see \cite{Ka00} for a more detailed
description of these boundaries in the case of usual time homogeneous random
walks. Note that Theorem 3.19 allows one to extend identification of
the Poisson boundary with the ``natural'' boundaries from usual time
homogeneous random walks to RWRTP for several other classes of groups,
including the polycyclic groups and the groups with infinitely many ends, cf.
\cite{Ka91} and \cite{Ka00}.

\head Appendix. Borel and Lebesgue spaces, conditional measures, discrete
equivalence relations and ergodic decompositions
\endhead

In this Appendix we sum up several fundamental definitions and facts about the
objects listed in its title, which are heavily used throughout the paper.
Although this language has become standard in the ergodic theory, it may be
less known to probabilists, which is why we preferred to expand on this rather
than to restrict ourselves just to an assortment of references. This brief
digest is based on the works of Rokhlin~\cite{Ro49}, Feldman--Moore \cite{FM77}
and Greschonig--Schmidt \cite{GS01}, also see the summaries of the theory of
Lebesgue spaces in \cite{Ro67} and \cite{CFS82}.

The only ``novelty'' here is Theorem A.7 on the transitivity of the ergodic
decomposition for discrete equivalence relations (which is used in the article
for proving Proposition~1.11 on the decomposition of the tail boundary of the
global chain into the tail boundaries of conditional random chains). We could
not find a precise reference to this fact in the available literature, although
it is no doubt known to specialists (as a part of ``folklore'', very much
present in this area, cf. the introduction to \cite{GS01}).

\bigskip

Let $(X,\B)$ be a {\it standard Borel space}, i.e., a complete separable metric
space equipped with the $\s$-algebra $\B$ of Borel subsets. It is isomorphic,
as a measurable space, to the unit interval with its Borel $\s$-algebra. The
measure space $(X,\F,m)$, where $m$ is a probability measure on $(X,\B)$, and
$\F$ is the {\it $m$-completion} of the $\s$-algebra $\B$ (i.e., $\F$ is
generated by $\B$ and all subsets in $X$ of $m$-measure $0$), is called a {\it
Lebesgue space}. Any purely non-atomic Lebesgue space is isomorphic to the unit
interval with the Lebesgue measure on it. The notion of a Lebesgue space can
also be defined intrinsically in the measure-theoretical category, without
referring to standard Borel spaces, see the fundamental paper of Rokhlin
\cite{Ro49} (for the lack of a more modern comprehensive exposition). As we are
working in the measure category, we shall always deal with Lebesgue spaces {\it
mod 0\/}, i.e., up to subsets of measure 0.

As usually, a {\it homomorphism} $\pi:(X,\F,m)\to (X',\F',m')$ of Lebesgue
spaces is a measure preserving surjection, i.e., for any subset $A'\in\F'$ its
preimage $\pi^{-1}(A')$ belongs to $\F$, and $m(\pi^{-1}(A'))=m'(A')$
(actually, any homomorphic image of a Lebesgue space is also a Lebesgue space).
Denote by $\F_\pi\subset\F$ the completed preimage $\s$-algebra
$\pi^{-1}(\F')$. An important property of Lebesgue spaces is that, conversely,
any complete sub-$\s$-algebra of $\F$ can be obtained in this way:

\proclaim{Theorem A.1 {\rm (Homomorphisms and sub-$\s$-algebras, \cite{Ro49})
}} Let $(X,\F,m)$ be a \linebreak Lebesgue space. Then for any complete
sub-$\s$-algebra $\F_1\subset\F$ there exists a unique (mod 0) homomorphism
$\pi:(X,\F,m)\to (X',\F',m')$ such that $\F_1=\F_\pi$. We shall call the space
$(X',\F',m')$ the {\sss quotient} of $(X,\F,m)$ by the $\s$-algebra $\F_1$.
\endproclaim

There is yet another class of objects in one-to-one correspondence with
complete sub-$\s$-algebras in Lebesgue spaces, that of {\it measurable
partitions}, which consists of the partitions $\xi_\pi$ of the space $(X,\F,m)$
into the preimages of a homomorphism of Lebesgue spaces $\pi:(X,\F,m)\to
(X',\F',m')$. Once again, there exists an intrinsic definition of measurable
partitions, see \cite{Ro49}; any such partition $\xi$ uniquely determines a
homomorphism $\pi$ of Lebesgue spaces with $\xi=\xi_\pi$ and, therefore, the
corresponding complete sub-$\s$-subalgebra.

Note that uniqueness in Theorem A.1 implies that replacing the measure $m$ on
the space $(X,\F)$ with an equivalent measure $\wt m$ leads to the same
quotient space $(X',\F')$ with the measure $\wt m'$ equivalent to $m'$. This
observation allows one to talk in the setup of Theorem A.1 about quotients for
spaces with $\s$-finite measures as well.

\bigskip

For any measurable subset $C$ {\it (condition)} in a probability space
$(X,\F,m)$ with $m(C)>0$ one has the corresponding {\it conditional measure}
$m_C(A)=m(A)/m(C)$ defined on the $\s$-algebra $\F_C=\{A\cap C: A\in\F\}$ of
subsets of $C$. In Lebesgue spaces the notion of a conditional measure can be
extended to measure 0 conditions as well, namely, one can define in a coherent
way conditional measures with respect to almost all elements of any measurable
partition.

\proclaim{Theorem A.2 {\rm (Canonical systems of conditional measures,
\cite{Ro49})}} Given a homomorphism $\pi:(X,\F,m)\to(X',\F',m')$ of Lebesgue
spaces, denote the elements of the preimage partition $\xi_\pi$ by
$X_{x'}=\pi^{-1}(x'),\,x'\in X'$ . Then the sets $X_{x'}$ can be endowed with
structures of Lebesgue spaces $(X_{x'},\F_{x'},m_{x'})$ in such a way that for
any subset $A\in\F$ the intersections $A_{x'}=A\cap X_{x'}$ belong to $\F_{x'}$
for $m'$-a.e. $x'\in X'$, the function $x'\mapsto m_{x'}(A_{x'})$ is measurable
with respect to the $\s$-algebra $\F'$, and
$$
m(A) = \int m_{x'}(A_{x'})\,dm'(x') \;.
$$
We shall write the above identity as the decomposition of the measure $m$ into
an integral of the conditional measures $m_{x'}=m_{X_{x'}}$:
$$
m = \int m_{x'}\, dm'(x') \;.
$$
The system of measure spaces $(X_{x'},\F_{x'},m_{x'})$ with the above property
is (mod 0) unique, and it is called the {\sss canonical system of conditional
measures} of the measure $m$ with respect to the homomorphism $\pi$ (or: with
respect to the measurable partition $\xi_\pi$, with respect to the complete
sub-$\s$-algebra $\F_\pi$).
\endproclaim

Let three Lebesgue spaces $(X,\F,m), (X',\F',m')$ and $(X'',\F'',m'')$ be
connected with homomorphisms $\pi',\pi''$ and $p$ as shown in the following
commutative diagram:
$$
\beginpicture
 \setcoordinatesystem units <1mm,1mm>
 \setplotarea x from -4 to 50, y from -20 to 0
 \put{$(X,\F,m)$} at 0 0
 \put{$(X',\F',m')$} at -1 -16
 \put{$(X'',\F'',m'')\qquad,$} at 41 -16
 \put{$\pi'$} at -4 -7
 \put{$\pi''$} at 21 -7
 \put{$p$} at 16 -14
 \arrow <3.6truept> [.4,1.75] from -1 -4 to -1 -12
 \arrow <3.6truept> [.4,1.75] from 10 -4 to 24 -12
 \arrow <3.6truept> [.4,1.75] from 10 -16 to 24 -16
\endpicture
$$
{\vskip -2cm \leftline{(A.1)} \vskip 1.5cm \parindent 0pt so that the partition
$\xi_{\pi''}$ is {\it smaller} than the partition $\xi_{\pi'}$ (its elements
are bigger). Then for a.e. element $X_{x''},\,x''\in X''$ of the partition
$\xi_{\pi''}$ the partition of the space $(X_{x''},\F_{x''},m_{x''})$ induced
by $\xi_{\pi'}$ is measurable (its elements are the elements $X_{x'}$ of
$\xi_{\pi'}$ contained in $X_{x''}$, i.e., those, for which $p(x')=x''$) and is
determined by the restriction of the homomorphism $\pi'$ to $X_{x''}$. Thus,
any element $X_{x'}$ of the partition $\xi_{\pi'}$ can also be considered as an
element of the induced partition of the space $X_{x''},\, x''=p(x')$. Then the
uniqueness part of Theorem A.2 implies:}

\proclaim{Theorem A.3 {\rm (Transitivity of canonical systems of conditional
measures, \cite{Ro67, \S 1.7})}} \linebreak For $m'$-a.e. $x'\in X'$ the
Lebesgue space structure on the set $X_{x'}$ induced by conditioning directly
by the partition $\xi_{\pi'}$ and the structure induced conditioning first by
the smaller partition $\xi_{\pi''}$ and then by the induced partition of
$X_{x''}$ are the same, i.e.,
$$
\bigl(m_{p(x')}\bigr)_{x'} = m_{x'} \qquad \text{\rm for $m'$-a.e.\;} x'\in X'
\;.
$$
Therefore, for $m''$-a.e. $x''\in X''$
$$
m_{x''} = \int m_{x'} \, d m'_{x''} (x') \;, \tag A.2
$$
where $m'_{x''}$ are the conditional measures of $m'$ determined by the
projection $p$.
\endproclaim

An {\it equivalence relation} $R\subset X\times X$ on a standard Borel space
$(X,\B)$ is called {\it Borel} if $R$ is a Borel subset of the space $X\times
X$ (endowed with the product Borel structure). It is called {\it countable} if
all its classes are at most countable. If $G$ is a countable group of Borel
automorphisms of $(X,\B)$, then the associated {\it orbit equivalence relation}
$R_G=\{(x,gx):x\in X,\,g\in G\}$ is obviously countable and Borel. The converse
is also true:

\proclaim{Theorem A.4 {\rm (General and orbit equivalence relations,
\cite{FM77})}} For any countable Borel equivalence relation $R$ on a standard
Borel space $(X,\B)$ there exists a countable group $G$ of Borel automorphisms
of $(X,\B)$ such that $R=R_G$.
\endproclaim

A measure $m$ on $(X,\B)$ is called {\it $R$-quasi-invariant} (and the
equivalence relation $R$ is then called {\it non-singular} with respect to
$m$), if it is quasi-invariant with respect to a certain group $G$ with
$R=R_G$. This property does not depend on the choice of the group $G$ and can
actually be formulated for Lebesgue (rather than Borel) spaces in intrinsic
terms without evoking Theorem A.4. In the measure-theoretical category we shall
refer to (classes mod 0) of non-singular countable Borel equivalence relations
as {\it discrete equivalence relations}.

\remark{Remark A.5} By using Theorem A.4 it is easy to see that for any
countable Borel equivalence relation $R$ and a Borel measure $m$ there always
exists a Borel $R$-quasi-invariant measure $\wt m$ dominating $m$.
\endremark

Let $R$ be a discrete equivalence relation on a Lebesgue space $(X,\F,m)$. A
subset $A\in\F$ is called {\it $R$-invariant} if it is (mod 0) a union of
$R$-classes. Denote by $\F_R$ the complete sub-$\s$-algebra of measurable
$R$-invariant sets. If the $\s$-algebra $\F_R$ is trivial, i.e., if any
measurable $R$-invariant set is (mod 0) either trivial or coincides with $X$,
then the equivalence relation $(X,\F,m,R)$ is called {\it ergodic} (or: the
measure space $(X,\F,m)$ and the measure $m$ are called {\it $R$-ergodic}).

\proclaim{Theorem A.6 {\rm (Ergodic decomposition, see \cite{GS01} and the
references therein)}} Let $R$ be a discrete equivalence relation on a Lebesgue
space $(X,\F,m)$. Denote by $\pi_R:(X,\F,m)\to (X',\F',m')$ the homomorphism of
the space $(X,m,\F)$ determined by the $\s$-algebra $\F_R$ of $R$-invariant
sets. Let $\{(X_{x'},\F_{x'},m_{x'})\}$ be the associated canonical system of
conditional measures. Then for $m'$-a.e. $x'\in X'$ the set
$X_{x'}=\pi_R^{-1}(x') $ as a union of $R$-classes, and the restriction
$R_{x'}$ of $R$ onto $X_{x'}$ is an ergodic discrete equivalence relation on
the space $(X_{x'},\F_{x'},m_{x'})$. The homomorphism $\pi_R$ is the unique
homomorphism of the space $(X,\F,m)$ with these properties. The spaces
$(X_{x'},\F_{x'},m_{x'})$ are called the {\sss ergodic components} of the space
$(X,\F,m)$ with respect to the equivalence relation $R$, and the quotient space
$(X',\F',m')$ is called the {\sss space of ergodic components}.
\endproclaim

\proclaim{Theorem A.7 {\rm (Transitivity of the ergodic decomposition)}} Let
$R$ be a discrete equivalence relation on a Lebesgue space $(X,\F,m)$. Denote
by $(X',\F',m')=\pi'(X,\F,m)$ the associated space of ergodic components, and
let $(X'',\F'',m'')=p(X',\F',m')=\pi''(X,\F,m)$ be a quotient of the space
$(X',\F',m')$, see the diagram {\rm (A.1).} Denote by $m_{x'},\,x'\in X'$ and
$m_{x''},\,x''\in X''$ the conditional measures of $m$ determined by the
homomorphisms $\pi'$ and $\pi''$, respectively, and by $m'_{x''},\,x''\in X''$
the conditional measures of $m'$ determined by the homomorphism $\pi$. Then
formula {\rm (A.2)} provides the ergodic decompositions of the measures
$m_{x''}$ with respect to the corresponding restricted equivalence relations
$R_{x''}$.
\endproclaim

\demo{Proof} This is a direct consequence of the uniqueness of the ergodic
decomposition (Theorem A.6) and of the transitivity of the canonical systems of
conditional measures (Theorem A.3). Indeed, almost all measures $m_{x'}$ are
ergodic, so that formula (A.2) provides a decomposition of the measures
$m_{x''}$ as integrals of ergodic ones, which by Theorem A.6 is unique.
\demoend

\enddocument

\bye